\documentclass[12pt]{amsart}
\usepackage{amssymb}
\usepackage{amsfonts}
\usepackage{latexsym,amsmath,amsthm,amssymb,amsfonts}

\numberwithin{equation}{section}
\allowdisplaybreaks

\def\endproof{$\hfill\Box$\\}
\def\s{\,\,\,\,}

\numberwithin{equation}{section}
\newtheorem{theorem}{Theorem}[section]
\newtheorem{lem}[theorem]{Lemma}
\newtheorem{thm}[theorem]{Theorem}

\newtheorem{rem}[theorem]{Remark}

\newcounter{Cnumber}

\title[ ]
{\bf Landau-Lifshitz-Bloch equation on Riemannian manifold}
\author[ ]
{Boling Guo, \quad Zonglin Jia}

\date{}
\begin{document}
\maketitle

\begin{abstract}
In this article, we bring in Landau-Lifshitz-Bloch(LLB) equation on $m$-dimensional closed Riemannian manifold and prove that it admits a unique local solution. In addition, if $m\geqslant3$ and $L^{\infty}-$norm of initial data is sufficiently small, the solution can be extended globally. Moreover, if $m=2$, we can prove that the unique solution is global without assuming small initial data.
\end{abstract}

\section{Introduction}
Landau-Lifshitz-Gilbert equation describes physical properties of micromagnetic at temperatures below the critical temperature. The equation is as follows:
\begin{equation}\label{1.1}
\frac{\partial m}{\partial t}=\lambda_1m\times H_{eff}-\lambda_2m\times(m\times H_{eff})
\end{equation}
where $\times$ denotes the vector cross product in $\mathbb{R}^3$ and $H_{eff}$ is effective field while $\lambda_1$ and $\lambda_2$ are real constants.

However, at high temperature, the model must be replaced by following Landau-Lifshitz-Bloch equation(LLB)
\begin{equation}\label{1.2}
\frac{\partial u}{\partial t}=\gamma u\times H_{eff}+L_1\frac{1}{|u|^2}(u\cdot H_{eff})u-L_2\frac{1}{|u|^2}u\times(u\times H_{eff})
\end{equation}
where $\gamma$, $L_1$, $L_2$ are real numbers and $\gamma>0$. $H_{eff}$ is given by
\[H_{eff}=\Delta u-\frac{1}{\chi||}\Big(1+\frac{3T}{5(T-T_c)}|u|^2\Big)u.\]
where $T>T_c>0$ and $\chi||>0$.

Now let us recall some previous results about LLB. In \cite{L}, Le consider the case that $L_1=L_2=:\kappa_1>0$. At that time, he rewrites (\ref{1.2}) as
\begin{eqnarray}\label{1.31}
\frac{\partial u}{\partial t}=\kappa_1\Delta u+\gamma u\times\Delta u-\kappa_2(1+\mu|u|^2)u
\end{eqnarray}
with $\kappa_2:=\frac{\kappa_1}{\chi||}$ and $\mu:=\frac{3T}{5(T-T_c)}$ and assume that $\kappa_2$, $\gamma$, $\mu$ is positive. Le has proven that above equation with Neumann boundary value conditions has global weak solution(the weak solution here is different from ordinary one). Inspired by Le, in \cite{J} Jia introduces following equation
\begin{equation}\label{1.3}
\left\{
\begin{array}{llll}
\partial_t u=\kappa_1\Delta u+\gamma\nabla F(u)\times\Delta u-\kappa_2(1+\mu\cdot F(u))\nabla F(u)\s\s\mbox{in}\s\Omega\times(0,\infty)\\
\frac{\partial u}{\partial\nu}=0\s\s\mbox{on}\s\partial\Omega\times(0,\infty)\\
u(\cdot ,0)=u_0\s\s\mbox{in}\s\Omega
\end{array}
\right.
\end{equation}
where $\Omega$ is a regular bounded domain of $\mathbb{R}^d(d\leqslant3)$, $\nu$ is outer normal direction of $\partial\Omega$ and $F\in C^3(\mathbb{R}^3)$ is a known function. He calls it Generalized Landau-Lifshitz-Bloch equation(GLLB) and gets that (\ref{1.3}) admits a local strong solution provided $u_0\in W^{2,2}(\Omega,\mathbb{R}^3)$ and $\frac{\partial u_0}{\partial\nu}=0$. In \cite{GLZ}, Guo, Li and Zeng consider the coming LLB equation with initial condition
\begin{eqnarray} \label{eq:6}
\left\{ \begin{aligned}
         &u_t=\Delta u+u\times\Delta u-\lambda(1+\mu|u|^2)u\s\s\mbox{in}\s\mathbb{R}^d\times(0,T)\\
                  &u(,0)=u_0\s\s\mbox{in}\s\mathbb{R}^d,
                          \end{aligned} \right.
                          \end{eqnarray}
where the constant $\lambda,\mu>0$. They prove the existence of smooth solutions of (\ref{eq:6}) in $\mathbb{R}^2$ or $\mathbb{R}^3$. And a small initial value condition should be added in the latter case.\\

In this paper, we would like to introduce a equation similar with (\ref{eq:6}) on Riemannian manifold. Before getting to this, we should make some preparation.

Let $\pi:(E,h,D)\longrightarrow(M,g,\nabla)$ denote a smooth vector bundle over an $m-$dimensional smooth closed Riemannian manifold $(M,g,\nabla)$ with $rank(E)=3$. $g$ means Riemannian metric of $M$ and $\nabla$ is its Levi-Civita connection. $h$ and $D$ are respectively metric and connection of $E$ such that $Dh=0$. Sometimes we also write $h$ as $\langle\cdot,\cdot\rangle$.
\subsection{$k-$times continuously differentiable section}Suppose $\Gamma(E)$ is the set of all sections in $E$. Under arbitrary local frame $\{e_{\alpha}:1\leqslant\alpha\leqslant3\}$, a section $s\in\Gamma(E)$ can be written in the form of $s=s^{\alpha}\cdot e_{\alpha}$. If $s^{\alpha}$ is $k-$times continuously differentiable, then we say $s$ is $k-$times continuously differentiable. Since $E$ is smooth, $k-$times continuous differentiability is independent of the choice of local frame. Define
$$\Gamma^k(E):=\{s\in\Gamma(E):\mbox{$s$ is $k-$times continuously differentiable}\}.$$
\subsection{Orientable vector bundle}$E$ is called orientable if there exists an $\omega\in E^*\wedge E^*\wedge E^*$ such that $\omega$ is continuous and for all $p\in M$, $\omega(p)\not=0$, where $E^*$ is dual bundle of $E$.

Suppose $\{e_1,e_2,e_3\}$ is a frame of $E$. It is called adapted to the orientation $\omega$ if
$$\omega(e_1,e_2,e_3)>0.$$
From now on, we always assume that $E$ is orientable unless otherwise stated.
\subsection{Cross product on orientable vector bundle}Suppose $\omega$ is an orientation of $E$. $\{e_{\alpha}:1\leqslant\alpha\leqslant3\}$ is a local frame of $E$ which is adapted to $\omega$. For any $f_1,f_2\in\Gamma(E)$, we assume that $f_1:=f_1^{\alpha}\cdot e_{\alpha}$, $f_2:=f_2^{\alpha}\cdot e_{\alpha}$. Their cross product $\times$ is defined as follow
$$(f_1\times f_2)(p):=f_1(p)\times f_2(p),$$
where
\begin{eqnarray*}
f_1(p)\times f_2(p):&=&\big(f_1^2(p)\cdot f_2^3(p)-f_2^2(p)\cdot f_1^3(p)\big)\cdot e_1(p)\\
&+&\big(f_2^1(p)\cdot f_1^3(p)-f_1^1(p)\cdot f_2^3(p)\big)\cdot e_2(p)\\
&+&\big(f_1^1(p)\cdot f_2^2(p)-f_2^1(p)\cdot f_1^2(p)\big)\cdot e_3(p).
\end{eqnarray*}
It is not hard to verify that $f_1(p)\times f_2(p)$ does not depend upon the choice of local frames which are adapted to $\omega$.
\subsection{Laplace operator on vector bundle}
Define a functional $Energy$ on $\Gamma^{2}(E)$ which is given in the form of
$$Energy(X):=\frac{1}{2}\int_M|DX|^2\,dM.$$
It is not hard to see that the Euler-Lagrange equation of $Energy$ is
$$\Delta X:=g^{ij}\cdot(D^2X)\Big(\frac{\partial}{\partial x^i},\frac{\partial}{\partial x^j}\Big)=0,$$
where $g_{ij}:=g\Big(\frac{\partial}{\partial x^i},\frac{\partial}{\partial x^j}\Big)$ and $(g^{ij})$ is the inverse matrix of $(g_{ij})$. Then we say that $\Delta$ is the Laplace operator on vector bundle $E$.
\subsection{sections depending on time}
A section depending on time is a map
$$V:I\longrightarrow\Gamma(E),$$
where $I$ is an interval of $\mathbb{R}$. Under arbitrary local frame $\{e_{\alpha}:1\leqslant\alpha\leqslant3\}$, $V(t,x)$ can be written as $V(t,x):=V^{\alpha}(t,x)\cdot e_{\alpha}(x)$. If $V^{\alpha}$ is $k-$times continuously differentiable with respect to $t$, we say $V$ is $k-$times continuously differentiable with respect to $t$ and use the symbol $C^k(I,\Gamma(E))$ to denote all such $V$. Since $E$ is smooth, differentiability with respect to time is independent of the choice of local frame. Moreover, we define
$$(\partial_t^kV)(t,x):=(\partial_t^kV^{\alpha})(t,x)\cdot e_{\alpha}(x).$$
\subsection{Sobolev space on vector bundle} Equip $\Gamma^k(E)$ with a norm $||\cdot||_{H^{k,p}}(p\geqslant1)$ which is defined as follow
$$||s||^p_{H^{k,p}}:=\sum\limits_{i=0}^k\int_M|D^is|^p\,dM.$$
The Sobolev space $H^{k,p}(E)$ is the completion of $\Gamma^k(E)$ with respect to the norm $||\cdot||_{H^{k,p}}$. For convenience, we also denote $H^{k,2}$ by $H^k$ and $||\cdot||_{H^{0,p}}$ by $||\cdot||_p$.
\\

Having the above preparation, we will give the definition of Landau-Lifshitz-Bloch equation(LLB) on Riemannian manifold.

For any $T>0,\lambda>0$ and $\mu>0$, let us consider a section depending on time $V\in C^1([0,T],\Gamma^{2}(E))$. LLB is just the following equation
\begin{eqnarray}\label{eq:1}
\left\{ \begin{aligned}
         &\partial_tV=\Delta V+V\times\Delta V-\lambda\cdot(1+\mu\cdot|V|^2)V \s\s\s\s\mbox{in}\s(0,T]\times M\\
                  &V(0,\cdot)=V_0
                          \end{aligned} \right.
\end{eqnarray}
\\

Our main results are as follow:
\begin{thm}\label{thm4}
Let $\pi:(E,h,D)\longrightarrow(M,g,\nabla)$ denote a smooth vector bundle over an $m-$dimensional smooth closed Riemannian manifold $(M,g,\nabla)$ with $rank(E)=3$ and $Dh=0$. $E$ is orientable. Given $l\geqslant m_0+1$(Here $m_0:=[\frac{m}{2}]+3$ and $[q]$ is the integral part of $q$) and $V_0\in H^l(E)$, there is a $T^*=T^*(||V_0||_{H^{m_0}})>0$ and a unique solution $V$ of (\ref{eq:1}) satisfying that for any $0\leqslant j\leqslant[\frac{l}{\hat{m}}]$($\hat{m}:=\max\{2,[\frac{m}{2}]+1\}$) and $\alpha\leqslant l-\hat{m}j$,
\begin{eqnarray}\label{6}
\partial_t^jD^{\alpha}V\in L^{\infty}([0,T^*],L^2(E)).
\end{eqnarray}
Furthermore, if $V_0\in\Gamma^{\infty}(E)$, then $V\in C^{\infty}([0,T^*],\Gamma^{\infty}(E))$.
\end{thm}

\begin{thm}\label{thm5}
Let $\pi:(E,h,D)\longrightarrow(M,g,\nabla)$ denote a smooth vector bundle over an $m-$dimensional smooth closed Riemannian manifold $(M,g,\nabla)$ with $rank(E)=3$, $m\geqslant3$ and $Dh=0$. $E$ is orientable. For any $T>0$ and $N\geqslant m_0+1$, there exists a $\hat{B}_N>0$ such that for all $V_0\in H^N(E)$ with $||V_0||_{\infty}\leqslant\hat{B}_N$, there is a unique solution of (\ref{eq:1}) satisfying
\begin{eqnarray}\label{7}
\partial_t^jD^{\alpha}V\in L^{\infty}([0,T],L^2(E))\s\forall0\leqslant j\leqslant\Big[\frac{N}{\hat{m}}\Big]\s\forall\alpha\leqslant N-\hat{m}j
\end{eqnarray}
and
\begin{eqnarray}\label{9}
\partial_t^iD^{\beta}V\in L^{2}([0,T],L^2(E))\s\forall0\leqslant i\leqslant\Big[\frac{N+1}{\hat{m}+1}\Big]\s\forall\beta\leqslant N+1-(\hat{m}+1)i.
\end{eqnarray}
Furthermore, if $V_0\in\Gamma^{\infty}(E)$, then $V\in C^{\infty}([0,T],\Gamma^{\infty}(E))$.
\end{thm}

\begin{thm}\label{thm9}
Let $\pi:(E,h,D)\longrightarrow(M,g,\nabla)$ denote a smooth vector bundle over an $2-$dimensional smooth closed Riemannian manifold $(M,g,\nabla)$ with $rank(E)=3$ and $Dh=0$. $E$ is orientable. For any $T>0$,      $N\geqslant5$ and $V_0\in H^5(E)$, there is a unique solution of (\ref{eq:1}) satisfying
\begin{eqnarray*}
\partial_t^jD^{\alpha}V\in L^{\infty}([0,T],L^2(E))\s\forall0\leqslant j\leqslant\Big[\frac{N}{2}\Big]\s\forall\alpha\leqslant N-2j
\end{eqnarray*}
and
\begin{eqnarray*}
\partial_t^iD^{\beta}V\in L^{2}([0,T],L^2(E))\s\forall0\leqslant i\leqslant\Big[\frac{N+1}{3}\Big]\s\forall\beta\leqslant N+1-3i.
\end{eqnarray*}
Furthermore, if $V_0\in\Gamma^{\infty}(E)$, then $V\in C^{\infty}([0,T],\Gamma^{\infty}(E))$.
\end{thm}
\section{Notation and Preliminaries}
In the paper, we appoint that the same indices appearing twice means summing it. And $Q_1\lesssim Q_2$ implies there is a universal constant $C$ such that $Q_1\leqslant C\cdot Q_2$.
\subsection{Riemannian curvature tensor on vector bundle}Using the connection $D$ on $E$, we can define a tensor $R^E$ called Riemannian curvature tensor. For any $X,Y\in TM$ and $s\in\Gamma^2(E)$,
$$R^E(X,Y)s:=D_XD_Ys-D_YD_Xs-D_{[X,Y]}s.$$
Let $R^M$ be the Riemannian curvature tensor of $M$. Being going to represent $R^M$ and $R^E$ in local frame, we appoint $\frac{\partial}{\partial x^i}$ as $\partial_i$. Then,
$$R^M(\partial_i,\partial_j)\partial_r:=(R^M)^h_{ijr}\cdot\partial_h\s\s
\mbox{and}
\s\s R^E(\partial_i,\partial_j)e_{\beta}:=(R^E)^{\alpha}_{ij\beta}\cdot e_{\alpha}.$$
Now we give two tensors
$$\mathcal{R}^M:=(\mathcal{R}^M)_{ijkl}\cdot dx^i\otimes dx^j\otimes dx^k\otimes dx^l$$
and
$$\mathcal{R}^E:=(\mathcal{R}^E)^{\alpha\beta}_{ij}\cdot dx^i\otimes dx^j\otimes e_{\alpha}\otimes e_{\beta},$$
where
$$(\mathcal{R}^M)_{ijkl}:=(R^M)^h_{ijk}\cdot g_{hl}\s\s\s\s\mbox{and}\s\s\s\s(\mathcal{R}^E)^{\alpha\beta}_{ij}:=(R^E)^{\alpha}_{ij\theta}\cdot h^{\theta\beta}.$$
$(h_{\alpha\beta})$ is the metric matrix of $h$ and $(h^{\theta\beta})$ is its inverse matrix.
\subsection{Cross product of tensors}
We also want to introduce cross product between two tensors. Given $S\in\Gamma(T^*M^{\otimes k}\otimes E)$ and $T\in\Gamma(T^*M^{\otimes l}\otimes E)$, let us define
$$S\times T:=(S_{i_1\cdots i_k}\times T_{j_1\cdots j_l})\otimes dx^{i_1}\otimes\cdots\otimes dx^{i_k}\otimes dx^{j_1}\otimes\cdots\otimes dx^{j_l},$$
where
$$S_{i_1\cdots i_k}:=S(\partial_{i_1},\cdots,\partial_{i_k})\s\s\s\s \mbox{and}\s\s\s\s T_{j_1\cdots j_l}:=T(\partial_{j_1},\cdots,\partial_{j_l}).$$
It is easy to check
\begin{eqnarray}\label{15}
|S\times T|\leqslant|S|\cdot|T|
\end{eqnarray}
\subsection{Properties of cross product}
\begin{thm}\label{thm3}
For any $f_1,f_2\in\Gamma^1(E)$, we have
\begin{eqnarray}\label{1}
D(f_1\times f_2)=(Df_1)\times f_2+f_1\times(Df_2).
\end{eqnarray}
\end{thm}
\textbf{Proof.} Take any $p\in M$. Then there exists a neighbourhood $U$ and a positive number $\delta$ such that the following map
$$\exp_p: N_{\delta}\triangleq\{\hat{v}\in T_pM:||\hat{v}||<\delta\}\longrightarrow U$$
is a diffeomorphism. Take $v\in T_pM$ such that $||v||=1$. Define $\gamma_v(t):=\exp_p(tv)$, where $t\in[0,\delta)$. Now take arbitrary orthonormal basis $\{e_{p\alpha}:1\leqslant\alpha\leqslant3\}$ in $E_p$ which is adapted to $\omega$ and let it move parallelly along $\gamma_v$ to get $$\{e_{\alpha}(t,v):t\in[0,\delta),1\leqslant\alpha\leqslant3\}.$$ Clearly,
$$w(t):=\omega(e_1(t,v),e_2(t,v),e_3(t,v))>0,\s\s\s\s\forall t\in[0,\delta)$$
since $w$ is a continuous function with respect to $t$. In the next, let $v$ range all the direction in $T_pM$ to obtain $$\{e_{\alpha}(t,v):t\in[0,\delta),v\in T_pM,||v||=1,1\leqslant\alpha\leqslant3\}.$$
It is a orthonormal frame on $U$ which is adapted to $\omega$ and
\begin{eqnarray}\label{2}
(De_{\alpha})(p)=0.
\end{eqnarray}

Assume that $f_1=f_1^{\alpha}\cdot e_{\alpha}$ and $f_2=f_2^{\beta}\cdot e_{\beta}$. Then, (\ref{2}) yields
$$Df_1(p)=df_1^{\alpha}(p)\otimes e_{\alpha}(p)\s\s\s\s\mbox{and}\s\s\s\s Df_2(p)=df_2^{\beta}(p)\otimes e_{\beta}(p).$$
Recalling the definition of cross product, we have
\begin{eqnarray*}
f_1\times f_2:=\big(f_1^2\cdot f_2^3-f_2^2\cdot f_1^3\big)\cdot e_1+\big(f_2^1\cdot f_1^3-f_1^1\cdot f_2^3\big)\cdot e_2+\big(f_1^1\cdot f_2^2-f_2^1\cdot f_1^2\big)\cdot e_3.
\end{eqnarray*}
Therefore, since of (\ref{2}), one can get
\begin{eqnarray}\label{3}
[D(f_1\times f_2)](p):&=&\big[df_1^2(p)\cdot f_2^3(p)+f_1^2(p)\cdot df_2^3(p)-df_2^2(p)\cdot f_1^3(p)-f_2^2(p)\cdot df_1^3(p)\big]\otimes e_1(p)\nonumber\\
&+&\big[df_2^1(p)\cdot f_1^3(p)+f_2^1(p)\cdot df_1^3(p)-df_1^1(p)\cdot f_2^3(p)-f_1^1(p)\cdot df_2^3(p)\big]\otimes e_2(p)\\
&+&\big[df_1^1(p)\cdot f_2^2(p)+f_1^1(p)\cdot df_2^2(p)-df_2^1(p)\cdot f_1^2(p)-f_2^1(p)\cdot df_1^2(p)\big]\otimes e_3(p),\nonumber
\end{eqnarray}
\begin{eqnarray}\label{4}
[f_1\times(Df_2)](p)&=&f_1(p)\times(Df_2)(p)\nonumber\\
&=&\big[f_1^2(p)\cdot df_2^3(p)-df_2^2(p)\cdot f_1^3(p)\big]\otimes e_1(p)\nonumber\\
&+&\big[df_2^1(p)\cdot f_1^3(p)-f_1^1(p)\cdot df_2^3(p)\big]\otimes e_2(p)\\
&+&\big[f_1^1(p)\cdot df_2^2(p)-df_2^1(p)\cdot f_1^2(p)\big]\otimes e_3(p),\nonumber
\end{eqnarray}
and
\begin{eqnarray}\label{5}
[(Df_1)\times f_2](p)&=&(Df_1)(p)\times f_2(p)\nonumber\\
&=&\big[df_1^2(p)\cdot f_2^3(p)-f_2^2(p)\cdot df_1^3(p)\big]\otimes e_1(p)\nonumber\\
&+&\big[f_2^1(p)\cdot df_1^3(p)-df_1^1(p)\cdot f_2^3(p)\big]\otimes e_2(p)\\
&+&\big[df_1^1(p)\cdot f_2^2(p)-f_2^1(p)\cdot df_1^2(p)\big]\otimes e_3(p).\nonumber
\end{eqnarray}
This theorem follows easily from combining (\ref{3}) with (\ref{4}) and (\ref{5}).
\endproof

Because of (\ref{1}), it is easy to verify that
\begin{eqnarray}\label{14}
D(S\times T)=(DS)\times T+S\times(DT),
\end{eqnarray}
provided $S\in\Gamma^1(T^*M^{\otimes k}\otimes E)$ and $T\in\Gamma^1(T^*M^{\otimes l}\otimes E)$.
\subsection{Hamilton's notation}Suppose $k,l,p,q\in\mathbb{N}$, $S\in T^*M^{\otimes k}\otimes E^{\otimes p}$ and $T\in T^*M^{\otimes l}\otimes E^{\otimes q}$, where
$$E^{\otimes p}:=\underbrace{E\otimes\cdots\otimes E}\limits_{\mbox{$p-$times}}.$$
we will write $S\ast T$, following Hamilton \cite{Hamiltion}, to denote a tensor formed by contraction on some indices of $S\otimes T$ using the coefficients $g^{ij}$ or $h_{\alpha\beta}$.
\begin{thm}\label{thm1}
$$|S\ast T|\leqslant|S|\cdot|T|$$
\end{thm}
\textbf{Proof.} We will get the above formula in an orthonormal basis of $M$ and an orthonormal basis of $E$.
\begin{eqnarray*}
|S\ast T|^2&=&\sum\limits_{\substack{\mbox{free}\\ \mbox{indices}}}\Bigg(\sum\limits_{\substack{\mbox{contracted}\\ \mbox{indices}}}S^{\alpha_1\cdots\alpha_p}_{i_1\cdots i_k}\cdot T^{\beta_1\cdots\beta_q}_{j_1\cdots j_l}\Bigg)^2\\
&\leqslant&\sum\limits_{\substack{\mbox{free}\\ \mbox{indices}}}\Bigg[\sum\limits_{\substack{\mbox{contracted}\\ \mbox{indices}}}\Big(S^{\alpha_1\cdots\alpha_p}_{i_1\cdots i_k}\Big)^2\Bigg]\cdot \Bigg[\sum\limits_{\substack{\mbox{contracted}\\ \mbox{indices}}}\Big(T^{\beta_1\cdots\beta_q}_{j_1\cdots j_l}\Big)^2\Bigg]\\
&\leqslant&\Bigg[\sum\limits_{\substack{\mbox{free}\\ \mbox{indices}}}\sum\limits_{\substack{\mbox{contracted}\\ \mbox{indices}}}\Big(S^{\alpha_1\cdots\alpha_p}_{i_1\cdots i_k}\Big)^2\Bigg]\cdot \Bigg[\sum\limits_{\substack{\mbox{free}\\ \mbox{indices}}}\sum\limits_{\substack{\mbox{contracted}\\ \mbox{indices}}}\Big(T^{\beta_1\cdots\beta_q}_{j_1\cdots j_l}\Big)^2\Bigg]\\
&=&|S|^2\cdot|T|^2
\end{eqnarray*}
\endproof
Because we do not specifically illustrate which indices are contracted, we have to appoint that
$$S_1\ast T_1-S_2\ast T_2:=S_1\ast T_1+S_2\ast T_2.$$

We will use the symbol $\mathfrak{q}_s(T_1,\cdots,T_r)$ for a polynomial in the tensors $T_1,\cdots,T_r$ and their iterated covariant derivatives with the $\ast$ product like
$$\mathfrak{q}_s(T_1,\cdots,T_r):=\sum\limits_{j_1+\cdots+j_r=s}c_{j_1\cdots j_r}\cdot D^{j_1}T_1\ast\cdots\ast D^{j_r}T_r,$$
where for $1\leqslant i\leqslant r$, $T_i\in\Gamma^{j_i}(T^*M^{\otimes t_i}\otimes E^{\otimes q_i})$ and $c_{j_1\cdots j_r}$ are some universal constants.
\subsection{Ricci identity}
Given $s\in\Gamma^2(T^*M^{\otimes k}\otimes E)$, it is obvious to see that $s$ can be written as follow
$$s:=s^{\alpha}_{i_1\cdots i_k}\cdot dx^{i_1}\otimes\cdots\otimes dx^{i_k}\otimes e_{\alpha}.$$
We denote $Ds$ in the form of components
$$Ds:=s^{\alpha}_{i_1\cdots i_k,p}\cdot dx^{i_1}\otimes\cdots\otimes dx^{i_k}\otimes dx^p\otimes e_{\alpha}.$$
At some time, we also employ the coming convention
\begin{eqnarray}\label{10}
Ds:=s_{i_1\cdots i_k,p}\cdot dx^{i_1}\otimes\cdots\otimes dx^{i_k}\otimes dx^p.
\end{eqnarray}
Thanks to the above agreement, Ricci identity is conveniently represented in the next theorem.
\begin{thm}\label{thm2}
\begin{eqnarray*}
&&s^{\alpha}_{i_1\cdots i_k,pq}-s^{\alpha}_{i_1\cdots i_k,qp}\\
&=&\sum\limits_{l=1}^ks^{\alpha}_{i_1\cdots i_{l-1}hi_{l+1}\cdots i_k}\cdot(R^M)^h_{pqi_l}-s^{\beta}_{i_1\cdots i_k}\cdot(R^E)^{\alpha}_{pq\beta}\\
&=&k\cdot s\ast\mathcal{R}^M+s\ast\mathcal{R}^E.
\end{eqnarray*}
\end{thm}
\textbf{Proof.}
The proof is straightforward if one takes normal coordinates. So we omit it.
\endproof

Given $V\in\Gamma^{k+1}(E)$ and $S\in\Gamma^{k+1}(T^*M\otimes E)$, by Theorem \ref{thm2} and induction, the following formulas are easy.\\
\textbf{Formula 1.}\\
There exist $a_{ij}\in\mathbb{Z}$ and $b_{rl}\in\mathbb{Z}$ such that
\begin{eqnarray}\label{11}
V_{,pi_1\cdots i_k}-V_{,i_1\cdots i_kp}&=&\sum\limits_{i+j=k-1}a_{ij}\cdot D^iV\ast D^j\mathcal{R}^E+\sum\limits_{r+l=k-2}b_{rl}\cdot D^{r+1}V\ast\nabla^l\mathcal{R}^M\nonumber\\
&=&\mathfrak{q}_{k-1}(V,\mathcal{R}^E)+\mathfrak{q}_{k-2}(DV,\mathcal{R}^M)
\end{eqnarray}
\textbf{Formula 2.}\\
There exist $a_{ij}\in\mathbb{Z}$ and $b_{rl}\in\mathbb{Z}$ such that
\begin{eqnarray}\label{12}
S_{p,qi_1\cdots i_k}-S_{p,i_1\cdots i_kq}&=&\sum\limits_{i+j=k-1}a_{ij}\cdot D^iS\ast D^j\mathcal{R}^E+\sum\limits_{r+l=k-1}b_{rl}\cdot D^{r}S\ast\nabla^l\mathcal{R}^M\nonumber\\
&=&\mathfrak{q}_{k-1}(S,\mathcal{R}^E)+\mathfrak{q}_{k-1}(S,\mathcal{R}^M)
\end{eqnarray}
\subsection{Interpolation for sections}
We shall prove Gagliardo-Nirenberg inequality of sections on vector bundle.
\begin{thm}\label{thm6}
$(M,g)$ is a $m-$dimensional smooth closed Riemannian manifold. $(E,h,D)$ is a smooth vector bundle over $M$ with $Dh=0$. $rank(E)$ may not be 3 and $E$ may not be orientable. Let $T$ be a smooth section of $E$. Given $s\in\mathbb{R}^+$ and $j\in\mathbb{Z}^+$, we will have
\begin{eqnarray}\label{29}
||D^jT||_{\frac{2s}{l}}\leqslant C(m,s,k,j)\cdot||D^kT||^{\frac{j}{k}}_{\frac{2s}{l+k-j}}\cdot||T||^{1-\frac{j}{k}}_{\frac{2s}{l-j}},
\end{eqnarray}
provided $k\in[j,\infty)\cap\mathbb{Z}$, $l\in[1,s]\cap[j,s+j+1-k]\cap\mathbb{Z}$.
\end{thm}
\textbf{Proof.} Apply induction for $j$.\\
\textbf{Step 1:} When $j=1$, (\ref{29}) is equivalent to
\begin{eqnarray}\label{30}
||DT||_{\frac{2s}{l}}\leqslant C(m,s,k)\cdot||D^kT||^{\frac{1}{k}}_{\frac{2s}{l+k-1}}\cdot||T||^{1-\frac{1}{k}}_{\frac{2s}{l-1}},
\end{eqnarray}
for all $l\in[1,s]\cap[1,s+2-k]\cap\mathbb{Z}$. In order to show (\ref{30}), we use induction for $k$.

When $k=1$, (\ref{30}) holds obviously.

When $k=2$, by 12.1 Theorem of \cite{Hamiltion} we know (\ref{30}) holds.

Assume that for $2\leqslant\hat{k}\leqslant k$, we obtain
\begin{eqnarray*}
||DT||_{\frac{2s}{l}}\leqslant C_1(m,s,\hat{k})\cdot||D^{\hat{k}}T||^{\frac{1}{\hat{k}}}_{\frac{2s}{l+\hat{k}-1}}\cdot||T||^{1-\frac{1}{\hat{k}}}_{\frac{2s}{l-1}},
\end{eqnarray*}
provided $l\in[1,s]\cap[1,s+2-\hat{k}]\cap\mathbb{Z}$.

When $\hat{k}=k+1$, pick any $l\in[1,s]\cap[1,s+2-(k+1)]\cap\mathbb{Z}$. Clearly,
$$l+1\in[1,s]\cap[1,s+2-k]\cap\mathbb{Z},$$
since $k\geqslant2$. Using induction hypothesis, we get
\begin{eqnarray}\label{31}
||D^2T||_{\frac{2s}{l+1}}\leqslant C_2(m,s,k)\cdot||D^k(DT)||^{\frac{1}{k}}_{\frac{2s}{l+k}}\cdot||DT||^{1-\frac{1}{k}}_{\frac{2s}{l}}.
\end{eqnarray}
Because $1\leqslant l\leqslant s+2-(k+1)<s$, using induction hypothesis for $k=2$ gives
\begin{eqnarray}\label{32}
||DT||_{\frac{2s}{l}}\leqslant C(m,s)\cdot||D^2T||^{\frac{1}{2}}_{\frac{2s}{l+1}}\cdot||T||^{\frac{1}{2}}_{\frac{2s}{l-1}}.
\end{eqnarray}
Combing (\ref{31}) with (\ref{32}) yields
$$||DT||_{\frac{2s}{l}}\leqslant C_3(m,s,k)\cdot||D^{k+1}T||^{\frac{1}{2k}}_{\frac{2s}{l+k}}\cdot||DT||^{\frac{1}{2}(1-\frac{1}{k})}_{\frac{2s}{l}}\cdot||T||^{\frac{1}{2}}_{\frac{2s}{l-1}},$$
which implies
$$||DT||_{\frac{2s}{l}}\leqslant C(m,s,k+1)\cdot||D^{k+1}T||^{\frac{1}{k+1}}_{\frac{2s}{l+k}}\cdot||T||^{1-\frac{1}{k+1}}_{\frac{2s}{l-1}}.$$
\textbf{Step 2:} Suppose that for all the indices not greater than $j$, (\ref{29}) is true. Now we consider $j+1$. At this moment, we take any $k\in[j+1,\infty)\cap\mathbb{Z}$ and any $l\in[1,s]\cap[j+1,s+j+2-k]\cap\mathbb{Z}$. It is easy to deduce that
$$k-1\in[j,\infty)\cap\mathbb{Z}\s\s\s\s\mbox{and}\s\s\s\s l\in[1,s]\cap[j,s+j+1-(k-1)]\cap\mathbb{Z}.$$

Using induction hypothesis leads to
\begin{eqnarray}\label{33}
||D^j(DT)||_{\frac{2s}{l}}\leqslant C_1(m,s,k,j)\cdot||D^{k-1}(DT)||^{\frac{j}{k-1}}_{\frac{2s}{l+k-1-j}}\cdot||DT||^{1-\frac{j}{k-1}}_{\frac{2s}{l-j}}.
\end{eqnarray}
Since $l-j\in[1,s]\cap[1,s+2-k]\cap\mathbb{Z}$, by Step 1 we have
\begin{eqnarray}\label{34}
||DT||_{\frac{2s}{l-j}}\leqslant C(m,s,k)\cdot||D^kT||^{\frac{1}{k}}_{\frac{2s}{l-j+k-1}}\cdot||T||^{1-\frac{1}{k}}_{\frac{2s}{l-j-1}}.
\end{eqnarray}
Combining (\ref{33}) with (\ref{34}) gives
\begin{eqnarray*}
||D^{j+1}T||_{\frac{2s}{l}}&\leqslant& C(m,s,k,j+1)\cdot||D^kT||^{\frac{j}{k-1}}_{\frac{2s}{l+k-1-j}}\cdot||D^kT||^{\frac{1}{k}(1-\frac{j}{k-1})}_{\frac{2s}{l+k-j-1}}\cdot||T||^{(1-\frac{1}{k})(1-\frac{j}{k-1})}_{\frac{2s}{l-j-1}}\\
&=&C(m,s,k,j+1)\cdot||D^kT||^{\frac{j+1}{k}}_{\frac{2s}{l+k-1-j}}\cdot||T||^{1-\frac{j+1}{k}}_{\frac{2s}{l-j-1}}.
\end{eqnarray*}
This completes the proof.\endproof
\begin{thm}\label{thm7}
$(M,g)$ is a $m-$dimensional smooth closed Riemannian manifold. $(E,h,D)$ is a smooth vector bundle over $M$ with $Dh=0$. $rank(E)$ may not be 3 and $E$ may not be orientable. Let $T$ be a smooth section of $E$. If $r,q\geqslant2$, then there is a universal constant $C=C(m,r,q,j,k)$ such that
\begin{eqnarray}\label{35}
||D^jT||_p\leqslant C\cdot||D^kT||_r^{\frac{j}{k}}\cdot||T||_q^{1-\frac{j}{k}},
\end{eqnarray}
provided
$$1\leqslant j\leqslant k\s\s\s\s\mbox{and}\s\s\s\s\frac{k}{p}=\frac{j}{r}+\frac{k-j}{q}.$$
\end{thm}
\textbf{Proof.} We consider 3 cases.\\
\textbf{Case 1:} When $2\leqslant r<q\leqslant\infty$, there exist $s$ and $l$ such that
$$q=\frac{2s}{l-j}\s\s\s\s\mbox{and}\s\s\s\s r=\frac{2s}{l+k-j}.$$
Since
$$\frac{k}{p}=\frac{j}{r}+\frac{k-j}{q},$$
we have $p=\frac{2s}{l}$. From Theorem \ref{thm6} it follows that
$$
||D^jT||_{\frac{2s}{l}}\leqslant C(m,s,k,j)\cdot||D^kT||^{\frac{j}{k}}_{\frac{2s}{l+k-j}}\cdot||T||^{1-\frac{j}{k}}_{\frac{2s}{l-j}},
$$
which means
$$
||D^jT||_p\leqslant C(m,r,q,j,k)\cdot||D^kT||_r^{\frac{j}{k}}\cdot||T||_q^{1-\frac{j}{k}}.
$$
\textbf{Case 2:} When $2\leqslant q<r\leqslant\infty$, the proof is similar.\\
\textbf{Case 3:} When $2\leqslant q=r$, clearly we have $p=q=r$. From 12.1 Theorem in \cite{Hamiltion} it follows that
$$
||DT||_p\leqslant C(m,p)\cdot||D^2T||^{\frac{1}{2}}_p\cdot||T||_p^{\frac{1}{2}},
$$
which implies
$$
||D^jT||_p\leqslant C(m,p)\cdot||D^{j+1}T||^{\frac{1}{2}}_p\cdot||D^{j-1}T||_p^{\frac{1}{2}}.
$$
Let $f(j):=||D^jT||_p$. It is easy to check that $f$ meets the condition of 12.5 Corollary in \cite{Hamiltion}. Then we conclude this theorem.\endproof
\section{Proof of Theorem \ref{thm4}}
Given any $T>0$, define an operator
$$P:C^1([0,T],\Gamma^2(E))\longrightarrow C([0,T],\Gamma(E)),$$
here
$$P(V):=\partial_tV-\Delta V-V\times\Delta V-\lambda(1+\mu|V|^2)V.$$
It is not difficult to check that the leading coefficient of the linearised operator of $P$ meets Legendre-Hadamard condition. By Main Theorem 1 in page 3 of \cite{Baker} we know (\ref{eq:1}) admits a unique local smooth solution $V$ provided $V_0\in\Gamma^{\infty}(E)$.

In the sequel, we would like to know the lower bound of maximal existence time $T_{\max}$ of the above smooth solution. Our strategy is to deduce a Gronwall inequality. That is to say, we shall control $\frac{d}{dt}||V(t)||^2_{H^l}$. Before getting to this, it is important to obtain an upper bound of $||V(t)||_{\infty}$.

Taking inner product with $|V|^{p-2}V(p>2)$ in (\ref{eq:1}), and integrating the result over $M$, we get
\begin{eqnarray*}
\int_M|V|^{p-2}\langle V,\partial_tV\rangle\,dM&=&\int_M|V|^{p-2}\langle V,\Delta V\rangle\,dM-\lambda\int_M(1+\mu|V|^2)|V|^p\,dM\\
&\leqslant&-\int_M|V|^{p-2}\cdot|DV|^2\,dM-(p-2)\int_M|V|^{p-4}\cdot|\langle V,DV\rangle|^2\,dM\leqslant0.
\end{eqnarray*}
The left hand side of the above inequality is $\frac{1}{p}\frac{d}{dt}||V(t)||_p^p$, so this inequality means
\begin{eqnarray*}
||V(t)||_p\leqslant||V_0||_p\s\s\s\s\forall t\in[0,T_{\max}).
\end{eqnarray*}
Taking the limit $p\rightarrow\infty$ leads to
\begin{eqnarray}\label{8}
||V(t)||_{\infty}\leqslant||V_0||_{\infty}\s\s\s\s\forall t\in[0,T_{\max}).
\end{eqnarray}

Given $k\geqslant1$, recalling our appointment (\ref{10}), we have the next identity
\begin{eqnarray*}
\frac{1}{2}\frac{d}{dt}\int_M|D^kV|^2\,dM&=&\int_Mg^{pq}g^{i_1j_1}\cdots g^{i_kj_k}\langle V_{,j_1\cdots j_k},V_{,pqi_1\cdots i_k}\rangle\,dM\\
&&+\int_Mg^{pq}g^{i_1j_1}\cdots g^{i_kj_k}\langle V_{,j_1\cdots j_k},(V\times V_{,p})_{,qi_1\cdots i_k}\rangle\,dM\\
&&-\lambda\cdot\mu\int_Mg^{i_1j_1}\cdots g^{i_kj_k}\langle V_{,j_1\cdots j_k},(|V|^2V)_{,i_1\cdots i_k}\rangle\,dM\\
&&-\lambda\int_M|D^kV|^2\,dM
\end{eqnarray*}
Applying (\ref{11}) and (\ref{12}) to exchange the order of derivatives yields
\begin{eqnarray}\label{13}
\frac{1}{2}\frac{d}{dt}\int_M|D^kV|^2\,dM&=&-\int_M|D^{k+1}V|^2\,dM+\int_MD^kV\ast\mathfrak{q}_k(V,\mathcal{R}^E)\,dM\\
&&+\int_MD^kV\ast\mathfrak{q}_{k-1}(DV,\mathcal{R}^M)\,dM\nonumber\\
&&-\int_Mg^{pq}g^{i_1j_1}\cdots g^{i_kj_k}\langle V_{,qj_1\cdots j_k},(V\times V_{,p})_{,i_1\cdots i_k}\rangle\,dM\nonumber\\
&&-\int_MD^kV\ast\mathfrak{q}_{k-1}(V\times DV,\mathcal{R}^M)\,dM\nonumber\\
&&-\int_MD^kV\ast\mathfrak{q}_{k-1}(V\times DV,\mathcal{R}^E)\,dM\nonumber\\
&&-\lambda\int_M|D^kV|^2\,dM-\lambda\mu\int_M|V|^2\cdot|D^kV|^2\,dM\nonumber\\
&&-\lambda\mu\sum\limits_{i+j=k-1}b_{ij}\cdot\int_MD^kV\ast D^{i+1}(|V|^2)\ast D^jV\,dM\nonumber\\
&&+\int_M\mathfrak{q}_k(V,\mathcal{R}^E)\ast D^{k-1}(V\times DV)\,dM\nonumber\\
&&+\int_M\mathfrak{q}_{k-1}(DV,\mathcal{R}^M)\ast D^{k-1}(V\times DV)\,dM,\nonumber
\end{eqnarray}
here $b_{ij}\in\mathbb{Z}^+$ are some universal constants. Note that
\begin{eqnarray*}
g^{pq}g^{i_1j_1}\cdots g^{i_kj_k}\langle V_{,qj_1\cdots j_k},(V\times V_{,p})_{,i_1\cdots i_k}\rangle=\sum\limits_{i+j=k-1}a_{ij}\cdot D^{k+1}V\ast(D^{i+1}V\times D^{j+1}V)
\end{eqnarray*}
where $a_{ij}\in\mathbb{Z}^+$ are some universal constants. Taking norms on the right hand side of (\ref{13}) leads to
\begin{eqnarray*}
\frac{1}{2}\frac{d}{dt}\int_M|D^kV|^2\,dM&\leqslant&-\int_M|D^{k+1}V|^2\,dM+\sum\limits_{i=0}^kC_i\int_M|D^kV|\cdot|D^iV|\,dM\\
&&+\sum\limits_{i+j=k-1}a_{ij}\int_M|D^{k+1}V|\cdot|D^{i+1}V|\cdot|D^{j+1}V|\,dM\\
&&+\sum\limits_{i=0}^{k-1}\bar{C}_i\int_M|D^kV|\cdot|D^i(V\times DV)|\,dM\\
&&+\lambda\mu\sum\limits_{i+j=k-1}b_{ij}\int_M|D^kV|\cdot|D^{i+1}(|V|^2)|\cdot|D^jV|\,dM\\
&&+\sum\limits_{i=0}^kC_i\int_M|D^iV|\cdot|D^{k-1}(V\times DV)|\,dM,
\end{eqnarray*}
where $C_i$ and $\bar{C}_i$ depend upon $\mathcal{R}^M$, $\mathcal{R}^E$ and their covariant differentiations. Applying (\ref{15}) and (\ref{14}) yields
\begin{eqnarray}\label{16}
\frac{1}{2}\frac{d}{dt}\int_M|D^kV|^2\,dM&\leqslant&-\int_M|D^{k+1}V|^2\,dM+\sum\limits_{i+j=k-1}a_{ij}\int_M|D^{k+1}V|\cdot|D^{i+1}V|\cdot|D^{j+1}V|\,dM\nonumber\\
&&+\tilde{C}_k\Big\{||D^kV||_2\cdot||V||_{H^k}+\sum\limits_{0\leqslant r+q\leqslant k-1}\int_M|D^kV|\cdot|D^rV|\cdot|D^{q+1}V|\,dM\nonumber\\
&&+\sum\limits_{r+q+j=k}\int_M|D^kV|\cdot|D^rV|\cdot|D^qV|\cdot|D^jV|\,dM\\
&&+\sum\limits_{i=0}^k\sum\limits_{r+q=k-1}\int_M|D^iV|\cdot|D^rV|\cdot|D^{q+1}V|\,dM\Big\},\nonumber
\end{eqnarray}
where $\tilde{C}_k$ depends upon $\mathcal{R}^M$, $\mathcal{R}^E$ and their covariant differentiations.
\begin{lem}\label{lem1}
There is a $C'_{m_0}>0$ depending on $\mathcal{R}^M$, $\mathcal{R}^E$ and their covariant differentiations such that, for any $t\in[0,T_{\max})$, we have
\begin{eqnarray*}
\frac{d}{dt}||V(t)||^2_{H^{m_0}}\leqslant C'_{m_0}\cdot(1+||V_0||^2_{H^{m_0}})\cdot\{||V(t)||^2_{H^{m_0}}+||V(t)||^4_{H^{m_0}}\}.
\end{eqnarray*}
\end{lem}
\textbf{Proof.} Given $1\leqslant k\leqslant m_0$, we consider
\begin{eqnarray*}
&&\sum\limits_{r+q+j=k}\int_M|D^kV|\cdot|D^rV|\cdot|D^qV|\cdot|D^jV|\,dM\\
&=&\sum\limits_{\substack{r+q+j=k\\ \max\{r,q,j\}=k}}\int_M|D^kV|\cdot|D^rV|\cdot|D^qV|\cdot|D^jV|\,dM\\
&&+\sum\limits_{\substack{r+q+j=k\\ \max\{r,q,j\}\leqslant k-1}}\int_M|D^kV|\cdot|D^rV|\cdot|D^qV|\cdot|D^jV|\,dM.
\end{eqnarray*}
Clearly,
\begin{eqnarray*}
\sum\limits_{\substack{r+q+j=k\\ \max\{r,q,j\}=k}}\int_M|D^kV|\cdot|D^rV|\cdot|D^qV|\cdot|D^jV|\,dM\lesssim||D^kV||_2^2\cdot||V||^2_{\infty}\leqslant||D^kV||_2^2\cdot||V_0||^2_{\infty}.
\end{eqnarray*}
And we want to derive the following
\begin{eqnarray*}
&&\sum\limits_{\substack{r+q+j=k\\ \max\{r,q,j\}\leqslant k-1}}\int_M|D^kV|\cdot|D^rV|\cdot|D^qV|\cdot|D^jV|\,dM\\
&\leqslant&\sum\limits_{\substack{r+q+j=k\\ \max\{r,q,j\}\leqslant k-1}}||D^kV||_2\cdot||D^rV||_{p_r}\cdot||D^qV||_{p_q}\cdot||D^jV||_{p_j}.
\end{eqnarray*}
where $p_r$, $p_q$ and $p_j$, belonging to $[1,\infty]$, will be determined later and satisfy
\begin{eqnarray}\label{20}
\frac{1}{p_r}+\frac{1}{p_q}+\frac{1}{p_j}=\frac{1}{2}.
\end{eqnarray}
And then we employ Theorem 2.1 due to \cite{DW1} to obtain
\begin{eqnarray*}
||D^rV||_{p_r}\lesssim||V||^{a_r}_{H^{m_0}}\cdot||V||_2^{1-a_r}\leqslant||V||_{H^{m_0}},
\end{eqnarray*}
\begin{eqnarray*}
||D^qV||_{p_q}\lesssim||V||^{a_q}_{H^{m_0}}\cdot||V||_2^{1-a_q}\leqslant||V||_{H^{m_0}},
\end{eqnarray*}
and
\begin{eqnarray*}
||D^jV||_{p_j}\lesssim||V||^{a_j}_{H^{m_0}}\cdot||V||_2^{1-a_j}\leqslant||V||_{H^{m_0}}.
\end{eqnarray*}
We hope $p_r$, $p_q$ and $p_j$ meet the next conditions:\\
\textbf{Condition 1.}
\begin{eqnarray*}
\frac{1}{p_r}=\frac{r}{m}+\frac{1}{2}-a_r\cdot\frac{m_0}{m}\s\s\s\s\mbox{with}\s\s\s\s a_r\in\Big[\frac{r}{m_0},1\Big),
\end{eqnarray*}
which is equivalent to
\begin{eqnarray}\label{17}
\frac{1}{p_r}\in\Big(\frac{r-m_0}{m}+\frac{1}{2},\frac{1}{2}\Big].
\end{eqnarray}
\textbf{Condition 2.}
\begin{eqnarray*}
\frac{1}{p_q}=\frac{q}{m}+\frac{1}{2}-a_q\cdot\frac{m_0}{m}\s\s\s\s\mbox{with}\s\s\s\s a_q\in\Big[\frac{q}{m_0},1\Big),
\end{eqnarray*}
which is equivalent to
\begin{eqnarray}\label{18}
\frac{1}{p_q}\in\Big(\frac{q-m_0}{m}+\frac{1}{2},\frac{1}{2}\Big].
\end{eqnarray}
\textbf{Condition 3.}
\begin{eqnarray*}
\frac{1}{p_j}=\frac{j}{m}+\frac{1}{2}-a_j\cdot\frac{m_0}{m}\s\s\s\s\mbox{with}\s\s\s\s a_j\in\Big[\frac{j}{m_0},1\Big),
\end{eqnarray*}
which is equivalent to
\begin{eqnarray}\label{19}
\frac{1}{p_j}\in\Big(\frac{j-m_0}{m}+\frac{1}{2},\frac{1}{2}\Big].
\end{eqnarray}
We claim there exist $p_r$, $p_q$ and $p_j$ which are in $[1,\infty]$ and satisfy (\ref{20}), (\ref{17}), (\ref{18}) and (\ref{19}). Obviously, that this claim holds is equivalent to
\begin{eqnarray}\label{21}
\Big(\frac{r-m_0}{m}+\frac{1}{2}\Big)+\Big(\frac{q-m_0}{m}+\frac{1}{2}\Big)+\Big(\frac{j-m_0}{m}+\frac{1}{2}\Big)<\frac{1}{2}\Longleftrightarrow k<3m_0-m.
\end{eqnarray}
Since $k\leqslant m_0$ and $m_0>\frac{m}{2}$, (\ref{21}) is true. In other words,
\begin{eqnarray*}
\sum\limits_{\substack{r+q+j=k\\ \max\{r,q,j\}\leqslant k-1}}\int_M|D^kV|\cdot|D^rV|\cdot|D^qV|\cdot|D^jV|\,dM\lesssim||D^kV||_2\cdot||V||^3_{H^{m_0}}
\end{eqnarray*}
In conclusion,
\begin{eqnarray}\label{23}
&&\sum\limits_{r+q+j=k}\int_M|D^kV|\cdot|D^rV|\cdot|D^qV|\cdot|D^jV|\,dM\\
&\lesssim&||D^kV||_2\cdot||V||^3_{H^{m_0}}+||D^kV||_2^2\cdot||V_0||^2_{\infty}.\nonumber
\end{eqnarray}
For the other terms of (\ref{16}), using the same methods, we get similar estimations:\\
\textbf{Estimation 1.}
\begin{eqnarray*}
&&\sum\limits_{i+j=k-1}a_{ij}\int_M|D^{k+1}V|\cdot|D^{i+1}V|\cdot|D^{j+1}V|\,dM\\
&\lesssim&||D^{k+1}V||_2\cdot||D^kV||_2\cdot||DV||_{\infty}+\sum\limits_{\substack{i+j=k-1\\ \max\{i,j\}\leqslant k-2}}a_{ij}\cdot||D^{k+1}V||_2\cdot||D^{i+1}V||_{p_i}\cdot||D^{j+1}V||_{p_j}\\
&\lesssim&||D^{k+1}V||_2\cdot||D^kV||_2\cdot||DV||_{\infty}+||D^{k+1}V||_2\cdot||V||^2_{H^{m_0}}\\
&\lesssim&||D^{k+1}V||_2\cdot||D^kV||_2\cdot||V||_{H^{m_0}}+||D^{k+1}V||_2\cdot||V||^2_{H^{m_0}},
\end{eqnarray*}
\textbf{Estimation 2.}
\begin{eqnarray*}
\sum\limits_{0\leqslant r+q\leqslant k-1}\int_M|D^kV|\cdot|D^rV|\cdot|D^{q+1}V|\,dM\lesssim||V_0||_{\infty}\cdot||D^kV||^2_2+||D^kV||_2\cdot||V||^2_{H^{m_0}},
\end{eqnarray*}
\textbf{Estimation 3.}
\begin{eqnarray*}
\sum\limits_{i=0}^k\sum\limits_{r+q=k-1}\int_M|D^iV|\cdot|D^rV|\cdot|D^{q+1}V|\,dM\lesssim||V||_{H^k}\cdot||D^kV||_2\cdot||V_0||_{\infty}+||V||_{H^k}\cdot||V||^2_{H^{m_0}}.
\end{eqnarray*}
Summing $k$ from $0$ to $m_0$ gives
\begin{eqnarray*}
\frac{1}{2}\frac{d}{dt}||V||^2_{H^{m_0}}&\leqslant&-||V||^2_{H^{m_0+1}}+L_{m_0}\cdot||V||_{H^{m_0+1}}\cdot||V||^2_{H^{m_0}}\\
&&+\tilde{L}_{m_0}\cdot\Big\{||V||^2_{H^{m_0}}+||V_0||_{\infty}\cdot||V||^2_{H^{m_0}}\\
&&+||V||^3_{H^{m_0}}+||V_0||^2_{\infty}\cdot||V||^2_{H^{m_0}}+||V||^4_{H^{m_0}}\Big\},
\end{eqnarray*}
where $L_{m_0}$ is universal and $\tilde{L}_{m_0}$ depends on $\mathcal{R}^M$, $\mathcal{R}^E$ and their covariant differentiations. Then the result follows easily from Young's inequality and Sobolev embedding
$$||V_0||_{\infty}\lesssim||V_0||_{H^{m_0}}.$$
This completes the proof.
\endproof

Consider an ODE
\begin{eqnarray}\label{eq:2}
\left\{ \begin{aligned}
        \frac{df}{dt}&=C'_{m_0}\cdot(1+||V_0||^2_{H^{m_0}})\cdot(f+f^2)\\
                  f(0)&=||V_0||^2_{H^{m_0}}.
                          \end{aligned} \right.
                          \end{eqnarray}
Solving the above equation to get an expression of $f$, we know that the maximal existence time of the solution $f(\cdot,||V_0||_{H^{m_0}})$ to (\ref{eq:2}) is not smaller than
$$T^*:=\frac{1}{C'_{m_0}\cdot(1+||V_0||^2_{H^{m_0}})}\log\Big(\frac{1+2||V_0||^2_{H^{m_0}}}{2||V_0||^2_{H^{m_0}}}\Big).$$
And $f(t,||V_0||_{H^{m_0}})$ is monotone increasing with respect to $t$. In other words, for all $t\in[0,T^*]$,
$$f(t,||V_0||_{H^{m_0}})\leqslant f(T^*,||V_0||_{H^{m_0}})=1+2||V_0||^2_{H^{m_0}}.$$
By comparison principle of ODE, we know that for any $t\in[0,\min\{T_{\max},T^*\})$,
$$||V(t)||_{H^{m_0}}\leqslant\sqrt{1+2||V_0||^2_{H^{m_0}}}.$$

In the sequel, we focus on the case that $k$ is sufficiently big.
\begin{lem}\label{lem2}
When $k\geqslant m_0+1$, there is a $Q_k>0$ depending on $\mathcal{R}^M$, $\mathcal{R}^E$ and their covariant differentiations such that, for any $t\in[0,T_{\max})$, we have
\begin{eqnarray}\label{22}
&&\frac{d}{dt}||V(t)||^2_{H^k}+||V(t)||^2_{H^{k+1}}\\
&\leqslant&Q_k\cdot\Big[||V(t)||^2_{H^k}\cdot||V(t)||^2_{H^{m_0}}+||V(t)||^4_{H^{k-1}}+||V(t)||^2_{H^k}\cdot(1+||V_0||_{\infty}+||V_0||^2_{\infty})\nonumber\\
&&+||V(t)||^2_{H^{k-1}}\cdot(||V(t)||^2_{H^{m_0}}+||V(t)||^4_{H^{m_0}})+||V(t)||^6_{H^{k-1}}\Big].\nonumber
\end{eqnarray}
\end{lem}
\textbf{Proof.} Firstly, let us calculate one term of (\ref{16}). Applying the same method of (\ref{23}), one can see easily that there are $p_i$ belonging to $[1,\infty]$ such that the following inequalities hold
\begin{eqnarray}\label{24}
&&\sum\limits_{i+j=k-1}a_{ij}\int_M|D^{k+1}V|\cdot|D^{i+1}V|\cdot|D^{j+1}V|\,dM\\
&\lesssim&||D^{k+1}V||_2\cdot||D^kV||_2\cdot||DV||_{\infty}+||D^{k+1}V||_2\cdot||D^{k-1}V||_2\cdot||D^2V||_{\infty}\nonumber\\
&&+\sum\limits_{\substack{i+j=k-1\\ \max\{i,j\}\leqslant k-3}}a_{ij}\cdot||D^{k+1}V||_2\cdot||D^{i+1}V||_{p_i}\cdot||D^{j+1}V||_{p_j}\nonumber\\
&\lesssim&||D^{k+1}V||_2\cdot||D^kV||_2\cdot||V||_{H^{m_0}}+||D^{k+1}V||_2\cdot||D^{k-1}V||_2\cdot||V||_{H^{m_0}}\nonumber\\
&&+||D^{k+1}V||_2\cdot||V||^2_{H^{k-1}}.\nonumber
\end{eqnarray}
By the same procedure, we get the next estimations:\\
\textbf{Estimation 4.}
\begin{eqnarray}\label{25}
&&\sum\limits_{0\leqslant r+q\leqslant k-1}\int_M|D^kV|\cdot|D^rV|\cdot|D^{q+1}V|\,dM\\
&\lesssim&||D^kV||_2^2\cdot||V_0||_{\infty}+||D^kV||_2\cdot||D^{k-1}V||_2\cdot||V||_{H^{m_0}}\nonumber\\
&&+||D^kV||_2\cdot||D^{k-2}V||_2\cdot||V||_{H^{m_0}}+||D^kV||_2\cdot||V||^2_{H^{k-1}}.\nonumber
\end{eqnarray}
\textbf{Estimation 5.}
\begin{eqnarray}\label{26}
&&\sum\limits_{r+q+j=k}\int_M|D^kV|\cdot|D^rV|\cdot|D^qV|\cdot|D^jV|\,dM\\
&\lesssim&\int_M|D^kV|^2\cdot|V|^2\,dM+\int_M|D^kV|\cdot|D^{k-1}V|\cdot|DV|\cdot|V|\,dM\nonumber\\
&&+\sum\limits_{\substack{r+q+j=k\\ \max\{r,q,j\}\leqslant k-2}}\int_M|D^kV|\cdot|D^rV|\cdot|D^qV|\cdot|D^jV|\,dM\nonumber\\
&\lesssim&||D^kV||_2^2\cdot||V||^2_{\infty}+||D^kV||_2\cdot||D^{k-1}V||_2\cdot||DV||_{\infty}\cdot||V||_{\infty}+||D^kV||_2\cdot||V||^3_{H^{k-1}}\nonumber\\
&\lesssim&||D^kV||_2^2\cdot||V_0||^2_{\infty}+||D^kV||_2\cdot||D^{k-1}V||_2\cdot||V||^2_{H^{m_0}}+||D^kV||_2\cdot||V||^3_{H^{k-1}}.\nonumber
\end{eqnarray}
\textbf{Estimation 6.}
\begin{eqnarray}\label{27}
&&\sum\limits_{i=0}^k\sum\limits_{r+q=k-1}\int_M|D^iV|\cdot|D^rV|\cdot|D^{q+1}V|\,dM\\
&\lesssim&||V||_{H^k}\cdot||V||_{\infty}\cdot||D^kV||_2+||V||_{H^k}\cdot||D^{k-1}V||_2\cdot||DV||_{\infty}\nonumber\\
&&+||V||_{H^k}\cdot||D^{k-2}V||_2\cdot||D^2V||_{\infty}+||V||_{H^k}\cdot||V||^2_{H^{k-1}}\nonumber\\
&\lesssim&||V||_{H^k}\cdot||V_0||_{\infty}\cdot||D^kV||_2+||V||_{H^k}\cdot||D^{k-1}V||_2\cdot||V||_{H^{m_0}}\nonumber\\
&&+||V||_{H^k}\cdot||D^{k-2}V||_2\cdot||V||_{H^{m_0}}+||V||_{H^k}\cdot||V||^2_{H^{k-1}}\nonumber
\end{eqnarray}

Substituting (\ref{24}), (\ref{25}), (\ref{26}) and (\ref{27}) into (\ref{16}) and then summing $k$ lead to
\begin{eqnarray*}
\frac{1}{2}\frac{d}{dt}||V||^2_{H^k}&\leqslant&-||V||^2_{H^{k+1}}+\hat{Q}_k\cdot(||V||_{H^{k+1}}\cdot||V||_{H^{k}}\cdot||V||_{H^{m_0}}+||V||_{H^{k+1}}\cdot||V||^2_{H^{k-1}})\\
&&+\tilde{Q}_k\cdot(||V||^2_{H^{k}}+||V||^2_{H^{k}}\cdot||V_0||_{\infty}+||V||_{H^{k}}\cdot||V||_{H^{k-1}}\cdot||V||_{H^{m_0}}\\
&&+||V||_{H^{k}}\cdot||V||^2_{H^{k-1}}+||V||^2_{H^{k}}\cdot||V_0||^2_{\infty}+||V||_{H^{k}}\cdot||V||_{H^{k-1}}\cdot||V||^2_{H^{m_0}}\\
&&+||V||_{H^{k}}\cdot||V||^3_{H^{k-1}}),
\end{eqnarray*}
where $\tilde{Q}_k>0$ depends on $\mathcal{R}^M$, $\mathcal{R}^E$ and their covariant differentiations. Using Young's inequality, we conclude this theorem.
\endproof

Note that (\ref{22}) is linear for $||V||^2_{H^k}$. It is now clear that inductively using (\ref{22}) one can show the existence of $N_k=N(||V_0||_{H^k},Q_k,Q_{k-1},\cdots,Q_{m_0+1})$ for any $k\geqslant m_0+1$ such that
$$||V(t)||_{H^k}\leqslant N_k\s\s\s\s\forall t\in[0,\min\{T_{\max},T^*\}),$$
which implies
$$T_{\max}\geqslant T^*.$$\\

Now we return to prove Theorem \ref{thm4}. Define
$$h(x):=\frac{1}{C'_{m_0}\cdot(1+x)}\log\Big(\frac{1+2x}{2x}\Big)$$
and we observe that it is a monotone decreasing function. Given $l\geqslant m_0+1$ and $V_0\in H^l(E)$, there are $V_{i0}\in\Gamma^{\infty}(E)$ such that as $i\rightarrow\infty$,
$$V_{i0}\longrightarrow V_0\s\s\s\s\mbox{strongly in $H^l(E)$}.$$
By the above discussion we know there exist $$T^*_i\geqslant h(||V_{i0}||^2_{H^{m_0}})>0\s\s\s\s\mbox{and}\s\s\s\s V_i\in C^{\infty}([0,T^*_i),\Gamma^{\infty}(E))$$ such that
\begin{eqnarray} \label{eq:3}
\left\{ \begin{aligned}
         &\partial_tV_i=\Delta V_i+V_i\times\Delta V_i-\lambda\cdot(1+\mu|V_i|^2)V_i \\
                 &V_i(0,\cdot)=V_{i0},
                          \end{aligned} \right.
                          \end{eqnarray}
here $T^*_i$ is the maximal existence time of $V_i$. Obviously, when $i$ is enough large,
$$||V_{i0}||^2_{H^{m_0}}\leqslant||V_{0}||^2_{H^{m_0}}+1\s\s\s\s\mbox{and}\s\s\s\s||V_{i0}||_{H^{l}}\leqslant||V_{0}||_{H^{l}}+1,$$
which imply
$$T^*_i\geqslant h(||V_{0}||^2_{H^{m_0}}+1):=2\tilde{\delta}>0$$
and
$$||V_i(t)||_{H^l}\leqslant N(||V_{0}||_{H^{l}}+1,Q_k,Q_{k-1},\cdots,Q_{m_0+1})\s\s\s\s\forall t\in[0,\tilde{\delta}].$$
Then $V_i$ is a bounded sequence in $L^{\infty}([0,\tilde{\delta}],H^l(E))$. It is not hard to verify that $\partial_tV_i$ is a bounded sequence in $L^{\infty}([0,\tilde{\delta}],L^2(E))$. So there exists a $V\in L^{\infty}([0,\tilde{\delta}],H^{l}(E))$ and a subsequence which is still denoted by $\{V_i\}$ such that
$$V_i\rightharpoonup V\s\s\s\s\mbox{weakly $*$ in $L^{\infty}([0,\tilde{\delta}],H^{l}(E))$}.$$
By Aubin-Lions lemma, one can find a subsequence still denoted by $\{V_i\}$ such that
$$V_i\longrightarrow V\s\s\s\s\mbox{strongly in $L^{\infty}([0,\tilde{\delta}],H^{l-1}(E))$}.$$
Because $l-1\geqslant m_0$, $H^{l-1}(E)$ can be embedded into $\Gamma^2(E)$. In other words, $V$ is a solution to (\ref{eq:1}). Using LLB to transform time derivatives into spatial derivatives gives that for all $0\leqslant j\leqslant\big[\frac{l}{\hat{m}}\big]$ and all $\alpha\leqslant l-\hat{m}j$, we have
\begin{eqnarray}\label{28}
\partial^j_tD^{\alpha}V\in L^{\infty}([0,\tilde{\delta}],L^2(E)).
\end{eqnarray}
\begin{rem}
The proof of (\ref{28}) is easy if one employs induction for $j$.
\end{rem}
At last, since $l\geqslant[\frac{m}{2}]+4$, by the same method of Theorem 3 in \cite{GLZ} it is not difficult to know that the solution of (\ref{eq:1}) with initial data $V_0\in H^l(E)$ is unique. This completes the proof.
\endproof
\section{Proof of Theorem \ref{thm5}}
Now we focus on global existence of LLB. Suppose that $V$ is the local smooth solution of (\ref{eq:1}). Our trick is to deduce a uniform estimation for $||V||_{H^k}$. To this goal, firstly we should get a linear Gronwall inequality.

By (\ref{16}) and H\"older inequality, we have
\begin{eqnarray}\label{37}
\frac{1}{2}\frac{d}{dt}||D^kV||_2^2&\leqslant&-||D^{k+1}V||_2^2+\sum\limits_{i+j=k-1}a_{ij}\cdot||D^{k+1}V||_2\cdot\Big|\Big||D^{i+1}V|\cdot|D^{j+1}V|\Big|\Big|_2\nonumber\\
&&+\tilde{C}_k\cdot\Big\{||V||^2_{H^k}+||D^kV||_2\sum\limits_{0\leqslant r+q\leqslant k-1}\Big|\Big||D^rV|\cdot|D^{q+1}V|\Big|\Big|_2\\
&&||D^kV||_2\sum\limits_{r+q+j=k}\Big|\Big||D^rV|\cdot|D^qV|\cdot|D^jV|\Big|\Big|_2\nonumber\\
&&+||V||_{H^k}\sum\limits_{r+q=k-1}\Big|\Big||D^rV|\cdot|D^{q+1}V|\Big|\Big|_2\Big\}.\nonumber
\end{eqnarray}

For the second term on the right hand side of (\ref{37}),
$$
\Big|\Big||D^{i+1}V|\cdot|D^{j+1}V|\Big|\Big|_2\leqslant||D^{i+1}V||_{\frac{2k+2}{i+1}}\cdot||D^{k-i}V||_{\frac{2k+2}{k-i}},
$$
since $i+j=k-1$. Theorem \ref{thm6} implies
$$
||D^{i+1}V||_{\frac{2k+2}{i+1}}\lesssim||D^{k+1}V||_2^{\frac{i+1}{k+1}}\cdot||V||_{\infty}^{\frac{k-i}{k+1}}
$$
and
$$
||D^{k-i}V||_{\frac{2k+2}{k-i}}\lesssim||D^{k+1}V||_2^{\frac{k-i}{k+1}}\cdot||V||_{\infty}^{\frac{i+1}{k+1}}.
$$
So
\begin{eqnarray}\label{38}
&&\sum\limits_{i+j=k-1}a_{ij}\cdot||D^{k+1}V||_2\cdot\Big|\Big||D^{i+1}V|\cdot|D^{j+1}V|\Big|\Big|_2\\
&\leqslant&B_k\cdot||V||_{\infty}\cdot||D^{k+1}V||_2^2\leqslant B_k\cdot||V_0||_{\infty}\cdot||D^{k+1}V||_2^2,\nonumber
\end{eqnarray}
where $B_k$ is a universal constant. By the same way, we will get
\begin{eqnarray}\label{39}
&&\sum\limits_{0\leqslant r+q\leqslant k-1}\Big|\Big||D^rV|\cdot|D^{q+1}V|\Big|\Big|_2\leqslant\sum\limits_{0\leqslant r+q\leqslant k-1}||D^rV||_{\frac{2r+2q+2}{r}}\cdot||D^{q+1}V||_{\frac{2r+2q+2}{q+1}}\\
&\lesssim&\sum\limits_{0\leqslant r+q\leqslant k-1}\Big[||D^{r+q+1}V||_2^{\frac{r}{r+q+1}}\cdot||V||^{\frac{q+1}{r+q+1}}_{\infty}\Big]\cdot\Big[||D^{r+q+1}V||_2^{\frac{q+1}{r+q+1}}\cdot||V||^{\frac{r}{r+q+1}}_{\infty}\Big]\nonumber\\
&=&\sum\limits_{0\leqslant r+q\leqslant k-1}||D^{r+q+1}V||_2\cdot||V||_{\infty}\lesssim||V||_{H^k}\cdot||V_0||_{\infty}\nonumber
\end{eqnarray}
and
\begin{eqnarray}\label{40}
\sum\limits_{r+q=k-1}\Big|\Big||D^rV|\cdot|D^{q+1}V|\Big|\Big|_2\lesssim||D^kV||_{2}\cdot||V_0||_{\infty}.
\end{eqnarray}
Moreover, Theorem \ref{thm6} yields
\begin{eqnarray}\label{41}
&&\sum\limits_{r+q+j=k}\Big|\Big||D^rV|\cdot|D^qV|\cdot|D^jV|\Big|\Big|_2\leqslant\sum\limits_{r+q+j=k}||D^rV||_{\frac{2k}{r}}\cdot||D^qV||_{\frac{2k}{q}}\cdot||D^jV||_{\frac{2k}{j}}\\
&\lesssim&\sum\limits_{r+q+j=k}\Big(||D^kV||^{\frac{r}{k}}_2\cdot||V||^{\frac{q+j}{k}}_{\infty}\Big)\cdot\Big(||D^kV||^{\frac{q}{k}}_2\cdot||V||^{\frac{r+j}{k}}_{\infty}\Big)\cdot\Big(||D^kV||^{\frac{j}{k}}_2\cdot||V||^{\frac{q+r}{k}}_{\infty}\Big)\nonumber\\
&\lesssim&||D^kV||_2\cdot||V||^2_{\infty}\leqslant||D^kV||_2\cdot||V_0||^2_{\infty}.\nonumber
\end{eqnarray}
Substituting (\ref{38}), (\ref{39}), (\ref{40}) and (\ref{41}) into (\ref{37}) leads to
\begin{eqnarray}\label{42}
&&\frac{1}{2}\frac{d}{dt}||D^kV||_2^2+(1-B_k\cdot||V_0||_{\infty})\cdot||D^{k+1}V||_2^2\\
&\leqslant&G_k\cdot\{||V||^2_{H^k}+||D^kV||_2\cdot||V||_{H^k}\cdot||V_0||_{\infty}+||D^kV||_2^2\cdot||V_0||^2_{\infty}\}\nonumber\\
&\leqslant&G_k\cdot(1+||V_0||_{\infty}+||V_0||^2_{\infty})\cdot||V||^2_{H^k},\nonumber
\end{eqnarray}
where $G_k$ depends upon $\mathcal{R}^M$, $\mathcal{R}^E$ and their covariant differentiations. In the sequel, using Gronwall inequality gives the following theorem.
\begin{thm}\label{thm8}
Given $N\in\mathbb{N}$, there exists an $\tilde{B}_N>0$ such that if $||V_0||_{\infty}\leqslant\tilde{B}_N$, we will obtain
\begin{eqnarray}\label{36}
||D^kV(t)||_2^2+\int^t_0||D^{k+1}V(s)||_2^2\,ds\leqslant C_k(||V_0||_{H^k},\tilde{B}_N,t),
\end{eqnarray}
provided $0\leqslant k\leqslant N$ and $t\in[0,T_{\max})$. Here $C_k(x,y,t)$ is monotone increasing with respect to $x$ and $t$.
\end{thm}
\textbf{Proof.} Employ induction for $N$.

In the case $N=0$, let $\tilde{B}_0:=1$. Taking inner product with $V$ in (\ref{eq:1}) and then integrating the result over $M$, we get
$$
\frac{1}{2}\frac{d}{dt}||V(t)||_2^2+||DV(t)||_2^2+\lambda\int_M(1+\mu|V(t)|^2)\cdot|V(t)|^2\,dM=0
$$
which is equivalent to
\begin{eqnarray}\label{*}
&&||V(t)||_2^2+2\int_0^t||DV(s)||_2^2\,ds\\
&&+2\lambda\int_0^tds\int_M(1+\mu|V(s)|^2)\cdot|V(s)|^2\,dM=||V_0||_2^2.\nonumber
\end{eqnarray}

Assume that for all the indices not larger than $N$, (\ref{36}) holds. Now we consider $N+1$.

Take $\tilde{B}_{N+1}:=\min\{\tilde{B}_N,\frac{1}{2B_{N+1}}\}$. If $||V_0||_{\infty}\leqslant\tilde{B}_{N+1}$, (\ref{42}) gives
\begin{eqnarray}\label{43}
&&\frac{1}{2}\frac{d}{dt}||D^{N+1}V||_2^2+\frac{1}{2}||D^{N+2}V||_2^2\\
&\leqslant& G_{N+1}\cdot(1+\tilde{B}_{N+1}+\tilde{B}^2_{N+1})\cdot\Big\{||D^{N+1}V||_2^2+\sum\limits_{k=0}^NC_k(||V_0||_{H^k},\tilde{B}_N,t)\Big\}.\nonumber
\end{eqnarray}
Then this theorem follows easily from Gronwall inequality. This completes the proof.\endproof

Now we return to prove Theorem \ref{thm5}. Given $T>0$ and $N\geqslant m_0+1=[\frac{m}{2}]+4$, we take any $V_0\in H^N(E)$ with $||V_0||_{\infty}\leqslant\frac{1}{2}\tilde{B}_N:=\hat{B}_N$(This $\tilde{B}_N$ is from Theorem \ref{thm8}). Then there are $V_{0i}\in\Gamma^{\infty}(E)$ converging to $V_0$ strongly in $H^N(E)$.

Suppose $V_i$ satisfies
\begin{eqnarray}\label{eq:4}
\left\{ \begin{aligned}
         &\partial_tV_i=\Delta V_i+V_i\times\Delta V_i-\lambda\cdot(1+\mu\cdot|V_i|^2)V_i \\
                  &V_i(0,\cdot)=V_{0i}
                          \end{aligned} \right.
\end{eqnarray}
and its maximal existence time is $T^*_i$. As $i$ is large enough, we have $$||V_{0i}||_{\infty}\leqslant2||V_0||_{\infty}\leqslant\tilde{B}_N,\s\s\s\s||V_{0i}||_{H^{m_0}}\leqslant2||V_0||_{H^{m_0}}$$
and $||V_{0i}||_{H^{N}}\leqslant2||V_0||_{H^{N}}$.

If $T^*_i<T$, then by Theorem \ref{thm8},
\begin{eqnarray*}
||V_i(t)||_{H^{m_0}}^2+\int^t_0||DV_i(s)||_{H^{m_0}}^2\,ds\leqslant C_{m_0}(||V_{0i}||_{H^{m_0}},\tilde{B}_N,t)\leqslant C_{m_0}(2||V_{0}||_{H^{m_0}},\tilde{B}_N,T),
\end{eqnarray*}
provided $t\in[0,T^*_i)$. Review that in the proof of Theorem \ref{thm4} we have defined a monotone decreasing function
$$
h(x):=\frac{1}{C'_{m_0}\cdot(1+x)}\log\Big(\frac{1+2x}{2x}\Big).
$$
So for arbitrary $t\in[0,T^*_i)$, it is obvious to see
$$h(||V_i(t)||^2_{H^{m_0}})\geqslant h\big(C_{m_0}(2||V_{0}||_{H^{m_0}},\tilde{B}_N,T)\big):=\delta_0>0.$$
Now we bring in a new system
\begin{eqnarray}\label{eq:5}
\left\{ \begin{aligned}
         &\partial_t\hat{V}_i=\Delta \hat{V}_i+\hat{V}_i\times\Delta \hat{V}_i-\lambda\cdot(1+\mu\cdot|\hat{V}_i|^2)\hat{V}_i \s\s\mbox{in}\s\s \Big(T^*_i-\frac{\delta_0}{2},\infty\Big)\times M\\
                  &\hat{V}_i\Big(T^*_i-\frac{\delta_0}{2},\cdot\Big)=V_i\Big(T^*_i-\frac{\delta_0}{2},\cdot\Big)
                          \end{aligned} \right.
\end{eqnarray}
The maximal existence time of $\hat{V}_i$ is not smaller than $h(||V_i(T^*_i-\frac{\delta_0}{2})||^2_{H^{m_0}})$ which is not smaller than $\delta_0$. By the uniqueness we know that for any $t\in[T^*_i-\frac{\delta_0}{2},T^*_i)$, $\hat{V}_i(t)=V_i(t)$. It means that $V_i$ can be extended to $[0,T^*_i+\frac{\delta_0}{2})$. Because $T^*_i$ is maximal, we get a contradiction. So $V_i\in C^{\infty}([0,T],\Gamma^{\infty}(E))$ and for all $t\in[0,T]$,
$$
||V_i(t)||_{H^N}^2+\int^t_0||DV_i(s)||_{H^N}^2\,ds\leqslant C_N(2||V_0||_{H^N},\tilde{B}_N,T).
$$

By the same method we prove local well-posedness one can know there is a
$$V\in L^{\infty}([0,T],H^N(E))\cap L^2([0,T],H^{N+1}(E))$$
such that $V_i$ converges to $V$ strongly in $L^{\infty}([0,T],H^{N-1}(E))$(in the sense of picking subsequence). It means $V$ is a solution of LLB.\\

At last, we claim (\ref{7}) and (\ref{9}) are true. Since (\ref{7}) is easy, we only prove (\ref{9}).\\
\textbf{Proof.} Employ induction for $i$.

When $i=0$, (\ref{9}) holds.

Suppose that for all the indices not bigger than $i$, (\ref{9}) is true.

Now we consider $i+1$. Choose any $\beta\in[0,N+1-(\hat{m}+1)(i+1)]\cap\mathbb{Z}$. Applying $\partial_t^iD^{\beta}$ to both sides of (\ref{eq:1}), we get
\begin{eqnarray*}
\partial_t^{i+1}D^{\beta}V=\partial_t^iD^{\beta}\Delta V-\partial_t^iD^{\beta}(V\times\Delta V)-\lambda\cdot\partial_t^iD^{\beta}\big[(1+\mu|V|^2)V\big],
\end{eqnarray*}
which implies
\begin{eqnarray}\label{44}
&&\int_0^{T}||\partial_t^{i+1}D^{\beta}V(s)||_2^2\,ds\\
&\lesssim&\int_0^{T}||\partial_t^{i}D^{\beta+2}V(s)||_2^2\,ds+\int_0^{T}||\partial_t^{i}D^{\beta}V(s)||_2^2\,ds\nonumber\\
&&+\sum\limits_{i',\beta'}\int_0^{T}\Big|\Big||\partial_t^{i'}D^{\beta'}V(s)|\cdot|\partial_t^{i-i'}D^{\beta-\beta'+2}V(s)|\Big|\Big|^2_2\,ds\nonumber\\
&&+\sum\limits_{\substack{i_1+i_2+i_3=i\\ \beta_1+\beta_2+\beta_3=\beta}}\int_0^{T}\Big|\Big||\partial_t^{i_1}D^{\beta_1}V(s)|\cdot|\partial_t^{i_2}D^{\beta_2}V(s)|\cdot|\partial_t^{i_3}D^{\beta_3}V(s)|\Big|\Big|^2_2\,ds.\nonumber
\end{eqnarray}
Because $\hat{m}:=\max\{2,[\frac{m}{2}]+1\}$,
$$i'\leqslant i\leqslant\Big[\frac{N+1}{\hat{m}+1}\Big]\leqslant\Big[\frac{N}{\hat{m}}\Big]\s\s\s\s\mbox{and}\s\s\s\s\beta'\leqslant\beta\leqslant N+1-(\hat{m}+1)(i+1)\leqslant N-\hat{m}\cdot i'-\hat{m},$$
(\ref{7}) yields
$$||\partial_t^{i'}D^{\beta'}V(t)||_{\infty}\leqslant||\partial_t^{i'}D^{\beta'}V(t)||_{H^{\hat{m}}}<\infty\s\s\s\s\forall t\in[0,T].$$
And since
$$\hat{m}(i-i')+(\beta-\beta'+2)\leqslant\hat{m}\cdot i+\beta+2\leqslant\hat{m}\cdot i+N+1-(\hat{m}+1)(i+1)+2=N-i-\hat{m}+2\leqslant N-i\leqslant N,$$
by (\ref{7}) we have
$$
\int_0^{T}||\partial_t^{i-i'}D^{\beta-\beta'+2}V(s)||_2^2\,ds\leqslant\sup\limits_{t\in[0,T]}\{||\partial_t^{i-i'}D^{\beta-\beta'+2}V(t)||^2_2\}\cdot T<\infty.
$$
So
\begin{eqnarray*}
&&\int_0^{T}\Big|\Big||\partial_t^{i'}D^{\beta'}V(s)|\cdot|\partial_t^{i-i'}D^{\beta-\beta'+2}V(s)|\Big|\Big|^2_2\,ds\\
&\leqslant&\sup\limits_{t\in[0,T]}\{||\partial_t^{i'}D^{\beta'}V(t)||^2_{\infty}\}\cdot\int_0^{T}||\partial_t^{i-i'}D^{\beta-\beta'+2}V(s)||_2^2\,ds<\infty
\end{eqnarray*}
For the other terms on the right hand side of (\ref{44}), using similar method we know all of them are strictly smaller than $\infty$.

This completes the proof.\endproof
\section{Proof of Theorem \ref{thm9}}
In this section, we need some formulas. Their proofs are tedious. So we only list the results.

\textbf{Formula 3.} Suppose that $V\in\Gamma^2(E)$. Then we will obtain
\begin{eqnarray*}
||\Delta V||_2^2&=&||D^2V||_2^2+2\int_M\langle DV,DV\ast\mathcal{R}^E\rangle\,dM\\
&&+\int_M\langle DV,V\ast D\mathcal{R}^E\rangle\,dM+\int_M\langle DV,DV\ast\mathcal{R}^M\rangle\,dM.
\end{eqnarray*}
\begin{rem}
Formula 3 easily implies
\begin{eqnarray}\label{5.4}
||D^2V||_2^2\leqslant||\Delta V||_2^2+\eta\cdot(||DV||_2^2+||V||_2^2),
\end{eqnarray}
where $\eta$ depends on $\mathcal{R}^M$, $\mathcal{R}^E$ and their covariant derivatives.
\end{rem}
\textbf{Formula 4.} Given $V\in\Gamma^3(E)$,
\begin{eqnarray*}
||\Delta DV||_2^2&=&||D^3V||_2^2+3\int_M\langle D^2V,D^2V\ast\mathcal{R}^M\rangle\,dM+2\int_M\langle D^2V,D^2V\ast\mathcal{R}^E\rangle\,dM\\
&&+\int_M\langle D^2V,DV\ast\nabla\mathcal{R}^M\rangle\,dM+\int_M\langle D^2V,DV\ast D\mathcal{R}^E\rangle\,dM.
\end{eqnarray*}
\begin{rem}
From Formula 4 it follows that
\begin{eqnarray}\label{45}
||D^3V||_2^2\leqslant||\Delta DV||_2^2+\eta_2\cdot(||D^2V||_2^2+||DV||_2^2),
\end{eqnarray}
where $\eta_2$ depends on $\mathcal{R}^M$, $\mathcal{R}^E$ and their covariant derivatives. Since by (\ref{11}) we have
\begin{eqnarray*}
\Delta DV=D\Delta V+\mathfrak{q}_1(V,\mathcal{R}^E)+\mathfrak{q}_0(DV,\mathcal{R}^M),
\end{eqnarray*}
integration by parts and H\"older's inequality yield
\begin{eqnarray}\label{47}
||\Delta DV||_2^2\leqslant||D\Delta V||_2^2+\eta_3\cdot(||D^2V||_2^2+||DV||_2^2+||V||_2^2),
\end{eqnarray}
where $\eta_3$ depends on $\mathcal{R}^M$, $\mathcal{R}^E$ and their covariant derivatives. Substituting (\ref{47}) into (\ref{45}) gives
\begin{eqnarray}\label{5.5}
||D^3V||_2^2\leqslant||D\Delta V||_2^2+(\eta_2+\eta_3)\cdot(||D^2V||_2^2+||DV||_2^2+||V||_2^2).
\end{eqnarray}
\end{rem}
\textbf{Formula 5.} If $V\in\Gamma^4(E)$, then
\begin{eqnarray*}
||D^4V||_2^2&=&||\Delta^2V||_2^2+\int_M\langle\mathfrak{q}_3(V,\mathcal{R}^E),D^3V\rangle\,dM+\int_M\langle\mathfrak{q}_2(DV,\mathcal{R}^M),D^3V\rangle\,dM\\
&&+\int_M\langle\mathfrak{q}_1(DV,\mathcal{R}^E),\mathfrak{q}_1(DV,\mathcal{R}^E)\rangle\,dM+\int_M\langle\mathfrak{q}_1(DV,\mathcal{R}^E),\mathfrak{q}_1(DV,\mathcal{R}^M)\rangle\,dM\\
&&+\int_M\langle\mathfrak{q}_1(DV,\mathcal{R}^M),\mathfrak{q}_1(DV,\mathcal{R}^M)\rangle\,dM.
\end{eqnarray*}
\begin{rem}
Formula 5, H\"older's inequality and Young's inequality lead to
\begin{eqnarray}\label{5.6}
||D^4V||_2^2&\leqslant&||\Delta^2V||_2^2+\eta_5\cdot(||D^3V||_2^2+||D^2V||_2^2+||DV||_2^2+||V||_2^2)\nonumber\\
&\leqslant&||\Delta^2V||_2^2+\eta_6\cdot(||D\Delta V||_2^2+||\Delta V||_2^2+||DV||_2^2+||V||_2^2),
\end{eqnarray}
where we have used (\ref{5.4}), (\ref{5.5}) and $\eta_5$, $\eta_6$ depend on $\mathcal{R}^M$, $\mathcal{R}^E$ and their covariant derivatives.
\end{rem}
\textbf{Formula 6.} Assume that $V\in\Gamma^4(E)$. Then we get
\begin{eqnarray}\label{48}
||D^2\Delta V||_2^2&=&||\Delta^2V||_2^2+\int_M\langle\mathfrak{q}_1(\Delta V,\mathcal{R}^E),D\Delta V\rangle\,dM+\int_M\langle\mathfrak{q}_0(D\Delta V,\mathcal{R}^M),D\Delta V\rangle\,dM\nonumber\\
&\lesssim&||\Delta^2V||_2^2+\eta_8\cdot(||D\Delta V||_2^2+||\Delta V||_2^2),
\end{eqnarray}
where we have applied H\"older's inequality, Young's inequality and $\eta_8$ depends on $\mathcal{R}^M$, $\mathcal{R}^E$ and their covariant derivatives.\\

Now let us go on to prove Theorem \ref{thm9}. Suppose that $V\in C^{\infty}\big([0,T^*),\Gamma^{\infty}(E)\big)$ is the unique local smooth solution of (\ref{eq:1}), where $T^*$ is its maximal existence time. First of all, we shall estimate $||DV(t)||_{\infty}$ for all $t\in[0,T^*)$. By Sobolev embedding it is easy to see that we only need to get a uniform upper bound of $||DV(t)||_{H^2}$(Note that in this section $m=2$). Combining (\ref{5.4}) and (\ref{5.5}) one can know that we only need to estimate
$$||DV(t)||_2^2+||\Delta V(t)||_2^2+||D\Delta V(t)||_2^2.$$

Using the same method of (2.2) in \cite{GLZ} we can get
\begin{eqnarray}\label{5.3}
||DV(t)||_2^2+\int_0^t||\Delta V(s)||_2^2\,ds\leqslant\lambda^2\cdot(1+\mu||V_0||^2_{\infty})^2\cdot||V_0||_2^2\cdot t+||DV_0||_2^2.
\end{eqnarray}
\\

For $||\Delta V(t)||_2^2$, our trick is to deduce a Gronwall's inequality. (\ref{eq:1}) yields
\begin{eqnarray}\label{49}
&&\frac{1}{2}\frac{d}{dt}||\Delta V(t)||_2^2+\int_M|D\Delta V(t)|^2\,dM+\lambda\cdot||\Delta V(t)||_2^2\\
&=&-\int_M[DV(t)\times\Delta V(t)]\ast D\Delta V(t)\,dM-\lambda\mu\cdot\int_M\langle\Delta[|V(t)|^2\cdot V(t)],\Delta V(t)\rangle\,dM\nonumber\\
&\leqslant&\int_M|DV(t)|\cdot|\Delta V(t)|\cdot|D\Delta V(t)|\,dM+C\cdot||V(t)||^2_{\infty}\cdot(||\Delta V(t)||^2_2+||DV(t)||_4^2)\nonumber\\
&\leqslant&||DV(t)||_4\cdot||\Delta V(t)||_4\cdot||D\Delta V(t)||_2+C\cdot||V(t)||^2_{\infty}\cdot(||\Delta V(t)||^2_2+||DV(t)||_4^2),\nonumber
\end{eqnarray}
here $C$ is a universal constant. Theorem 2.1 of \cite{DW1}, (\ref{5.4}) and (\ref{5.5}) give
\begin{eqnarray}\label{5.1}
||DV(t)||_4&\lesssim&||DV(t)||^{\frac{1}{4}}_{H^2}\cdot||DV(t)||_2^{\frac{3}{4}}\nonumber\\
&\leqslant&\eta_4\cdot(||D\Delta V(t)||_2^{\frac{1}{4}}+||\Delta V(t)||_2^{\frac{1}{4}}+||DV(t)||_2^{\frac{1}{4}}+||V(t)||_2^{\frac{1}{4}})\cdot||DV(t)||_2^{\frac{3}{4}}
\end{eqnarray}
and
\begin{eqnarray}\label{5.2}
||\Delta V(t)||_4&\lesssim&||\Delta V(t)||^{\frac{1}{2}}_{H^1}\cdot||\Delta V(t)||_2^{\frac{1}{2}}\nonumber\\
&\lesssim&(||D\Delta V(t)||_2^{\frac{1}{2}}+||\Delta V(t)||_2^{\frac{1}{2}})\cdot||\Delta V(t)||_2^{\frac{1}{2}},
\end{eqnarray}
where $\eta_4$ depends upon $\eta_2$, $\eta_3$ and $\eta$. Thus we derive
\begin{eqnarray*}
&&||DV(t)||_4\cdot||\Delta V(t)||_4\cdot||D\Delta V(t)||_2\\
&\lesssim&\eta_4\cdot||DV(t)||_2^{\frac{3}{4}}\cdot||\Delta V(t)||_2^{\frac{1}{2}}\cdot||D\Delta V(t)||_2^{\frac{7}{4}}\\
&&+\eta_4\cdot||DV(t)||_2^{\frac{3}{4}}\cdot||\Delta V(t)||_2^{\frac{3}{4}}\cdot||D\Delta V(t)||_2^{\frac{3}{2}}\\
&&+\eta_4\cdot||DV(t)||_2\cdot||\Delta V(t)||_2^{\frac{1}{2}}\cdot||D\Delta V(t)||_2^{\frac{3}{2}}\\
&&+\eta_4\cdot||DV(t)||_2^{\frac{3}{4}}\cdot||\Delta V(t)||_2\cdot||D\Delta V(t)||_2^{\frac{5}{4}}\\
&&+\eta_4\cdot||DV(t)||_2^{\frac{3}{4}}\cdot||\Delta V(t)||_2^{\frac{5}{4}}\cdot||D\Delta V(t)||_2\\
&&+\eta_4\cdot||DV(t)||_2\cdot||\Delta V(t)||_2\cdot||D\Delta V(t)||_2\\
&&+\eta_4\cdot||V(t)||_2^{\frac{1}{4}}\cdot||DV(t)||^{\frac{3}{4}}_2\cdot||\Delta V(t)||_2^{\frac{1}{2}}\cdot||D\Delta V(t)||_2^{\frac{3}{2}}\\
&&+\eta_4\cdot||V(t)||_2^{\frac{1}{4}}\cdot||DV(t)||^{\frac{3}{4}}_2\cdot||\Delta V(t)||_2\cdot||D\Delta V(t)||_2\\
&\leqslant&\eta_4\gamma_1\cdot(||\Delta V(t)||_2^{\frac{1}{2}}\cdot||D\Delta V(t)||_2^{\frac{7}{4}}+||\Delta V(t)||_2^{\frac{3}{4}}\cdot||D\Delta V(t)||_2^{\frac{3}{2}}\\
&&+||\Delta V(t)||_2^{\frac{1}{2}}\cdot||D\Delta V(t)||_2^{\frac{3}{2}}+||\Delta V(t)||_2\cdot||D\Delta V(t)||_2^{\frac{5}{4}}\\
&&+||\Delta V(t)||_2^{\frac{5}{4}}\cdot||D\Delta V(t)||_2+||\Delta V(t)||_2\cdot||D\Delta V(t)||_2),
\end{eqnarray*}
where $\gamma_1$ relies on $||V_0||_{\infty}, ||V_0||_2, ||DV_0||_2$ and $t$. From Young's inequality it follows that
\begin{eqnarray}\label{50}
&&||DV(t)||_4\cdot||\Delta V(t)||_4\cdot||D\Delta V(t)||_2\\
&\leqslant&\gamma_2\cdot(||\Delta V(t)||_2^4+1)+\frac{1}{4}\cdot||D\Delta V(t)||_2^2,\nonumber
\end{eqnarray}
where $\gamma_2$ is dependent of $\mathcal{R}^M$, $\mathcal{R}^E$ and their covariant derivatives, $||V_0||_{\infty}, ||V_0||_2, ||DV_0||_2$ and $t$.
By the same way, we have
\begin{eqnarray}\label{51}
&&||V(t)||^2_{\infty}\cdot(||\Delta V(t)||_2^2+||DV(t)||_4^2)\\
&\lesssim&||V(t)||^2_{\infty}\cdot||\Delta V(t)||_2^2\nonumber\\
&&+\eta_4^2\cdot||V_0||^2_{\infty}\cdot(||D\Delta V(t)||_2^{\frac{1}{2}}+||\Delta V(t)||_2^{\frac{1}{2}}+||DV(t)||_2^{\frac{1}{2}}+||V(t)||_2^{\frac{1}{2}})\cdot||DV(t)||_2^{\frac{3}{2}}\nonumber\\
&\leqslant&\gamma_3\cdot(||\Delta V(t)||_2^2+1)+\frac{1}{4}||D\Delta V(t)||_2^2,\nonumber
\end{eqnarray}
where $\gamma_3$ is dependent of $\mathcal{R}^M$, $\mathcal{R}^E$ and their covariant derivatives, $||V_0||_{\infty}, ||V_0||_2, ||DV_0||_2$ and $t$. Substituting (\ref{50}) and (\ref{51}) into (\ref{49}) and Young's inequality lead to
\begin{eqnarray*}
\frac{1}{2}\frac{d}{dt}||\Delta V(t)||_2^2+\frac{1}{2}||D\Delta V(t)||_2^2+\lambda||\Delta V(t)||_2^2\leqslant\gamma_4\cdot(||\Delta V(t)||_2^4+1),
\end{eqnarray*}
here $\gamma_4$ relies on $\gamma_2$ and $\gamma_3$. The generalized Gronwall's inequality says that if
$$\frac{df}{dt}\leqslant C\cdot f\cdot g+C,$$
then
$$f\leqslant C\cdot\exp\Big(\int_0^tg(s)\,ds\Big)+C.$$
So if we replace $f$ and $g$ by $||\Delta V(t)||_2^2$ and note that (\ref{5.3}) implies the boundedness of $\int_0^tg(s)\,ds$, then
\begin{eqnarray}\label{56}
||\Delta V(t)||_2^2\leqslant\gamma_5
\end{eqnarray}
which implies
\begin{eqnarray}\label{Eq:1}
\int_0^t||D\Delta V(s)||_2^2\,ds\leqslant2\gamma_4\cdot(\gamma_5^2+1)\cdot t,
\end{eqnarray}
where $\gamma_5$ is dependent of $\mathcal{R}^M$, $\mathcal{R}^E$ and their covariant derivatives, $||V_0||_{\infty}, ||V_0||_2, ||DV_0||_2$, $||\Delta V_0||_2$ and $t$.\\

In the sequel, we are going to estimate $||D\Delta V(t)||_2$ for all $t\in[0,T^*)$. (\ref{eq:1}) gives
\begin{eqnarray*}
\frac{1}{2}\frac{d}{dt}||D\Delta V(t)||_2^2&=&-||\Delta^2V(t)||_2^2-\int_M\langle\Delta^2V(t),\Delta[V(t)\times\Delta V(t)]\rangle\,dM\\
&&+\lambda\int_M\langle\Delta[(1+\mu|V(t)|^2)\cdot V(t)],\Delta^2 V(t)\rangle\,dM\\
&=&-||\Delta^2V(t)||_2^2-2\int_M\Delta^2V(t)\ast[DV(t)\times D\Delta V(t)]\,dM\\
&&+\lambda\int_M\langle\Delta[(1+\mu|V(t)|^2)\cdot V(t)],\Delta^2 V(t)\rangle\,dM\\
&=&-||\Delta^2V(t)||_2^2-2\int_M\Delta^2V(t)\ast[DV(t)\times D\Delta V(t)]\,dM\\
&&+\lambda\mu\int_M\langle\Delta\big[|V(t)|^2\cdot V(t)\big],\Delta^2 V(t)\rangle\,dM-\lambda||D\Delta V(t)||_2^2.
\end{eqnarray*}
On the other hand, H\"older inequality yields
\begin{eqnarray}\label{52}
&&\Big|\int_M\Delta^2V(t)\ast[DV(t)\times D\Delta V(t)]\,dM\Big|\\
&\leqslant&||DV(t)||_{\frac{16}{5}}\cdot||D\Delta V(t)||_{\frac{16}{3}}\cdot||\Delta^2V(t)||_2\nonumber
\end{eqnarray}
and
\begin{eqnarray*}
&&\Big|\int_M\langle\Delta\big[|V(t)|^2\cdot V(t)\big],\Delta^2 V(t)\rangle\,dM\Big|\\
&\lesssim&||V(t)||^2_{\infty}\cdot||\Delta^2V(t)||_2\cdot||\Delta V(t)||_2+||V(t)||_{\infty}\cdot||\Delta^2V(t)||_2\cdot||DV(t)||_4^2.
\end{eqnarray*}
By Sobolev Embedding, we have
\begin{eqnarray}\label{53}
||DV(t)||_{\frac{16}{5}}\lesssim||DV(t)||^{\frac{1}{8}}_{H^3}\cdot||DV(t)||^{\frac{7}{8}}_2.
\end{eqnarray}
Combining (\ref{5.4}), (\ref{5.5}) and (\ref{5.6}) we arrive at
\begin{eqnarray}\label{54}
||DV(t)||^{\frac{1}{8}}_{H^3}&\leqslant&\eta_7\cdot(||\Delta^2V(t)||_2^{\frac{1}{8}}+||D\Delta V(t)||_2^{\frac{1}{8}}\\
&&+||\Delta V(t)||_2^{\frac{1}{8}}+||DV(t)||_2^{\frac{1}{8}}+||V(t)||_2^{\frac{1}{8}}),\nonumber
\end{eqnarray}
where $\eta_7$ is dependent of $\mathcal{R}^M$, $\mathcal{R}^E$ and their covariant derivatives. Moreover,
\begin{eqnarray*}
||D\Delta V(t)||_{\frac{16}{3}}\lesssim||D\Delta V(t)||^{\frac{5}{8}}_{H^1}\cdot||D\Delta V(t)||_2^{\frac{3}{8}}
\end{eqnarray*}
and
\begin{eqnarray*}
||D\Delta V(t)||_{H^1}\lesssim||D^2\Delta V(t)||_2+||D\Delta V(t)||_2.
\end{eqnarray*}
By (\ref{48}) we are led to
\begin{eqnarray*}
||D\Delta V(t)||_{H^1}\lesssim||\Delta^2V(t)||_2+(\sqrt{\eta_8}+1)\cdot||D\Delta V(t)||_2+\sqrt{\eta_8}\cdot||\Delta V(t)||_2,
\end{eqnarray*}
which implies
\begin{eqnarray}\label{55}
||D\Delta V(t)||_{\frac{16}{3}}&\lesssim&\big[||\Delta^2V(t)||^{\frac{5}{8}}_2+(\eta_8^{\frac{5}{16}}+1)\cdot||D\Delta V(t)||_2^{\frac{5}{8}}\\
&&+\eta_8^{\frac{5}{16}}\cdot||\Delta V(t)||_2^{\frac{5}{8}}\big]\cdot||D\Delta V(t)||_2^{\frac{3}{8}}.\nonumber
\end{eqnarray}
Furthermore, substituting (\ref{53}), (\ref{54}) and (\ref{55}) into (\ref{52}) we arrive at
\begin{eqnarray}\label{5.7}
&&\Big|\int_M\Delta^2V(t)\ast[DV(t)\times D\Delta V(t)]\,dM\Big|\\
&\lesssim&\eta_7\cdot(||\Delta^2V(t)||_2^{\frac{1}{8}}+||D\Delta V(t)||_2^{\frac{1}{8}}+||\Delta V(t)||_2^{\frac{1}{8}}+||DV(t)||_2^{\frac{1}{8}}+||V(t)||_2^{\frac{1}{8}})\nonumber\\
&&\cdot\big[||\Delta^2V(t)||^{\frac{5}{8}}_2+(\eta_8^{\frac{5}{16}}+1)\cdot||D\Delta V(t)||_2^{\frac{5}{8}}+\eta_8^{\frac{5}{16}}\cdot||\Delta V(t)||_2^{\frac{5}{8}}\big]\nonumber\\
&&\cdot||D\Delta V(t)||_2^{\frac{3}{8}}\cdot||\Delta^2V(t)||_2\cdot||DV(t)||_2^{\frac{7}{8}}.\nonumber
\end{eqnarray}
Substituting the upper bounds of $||V(t)||_2$, $||DV(t)||_2$ and $||\Delta V(t)||_2$ into (\ref{5.7}) leads to
\begin{eqnarray*}
&&\Big|\int_M\Delta^2V(t)\ast[DV(t)\times D\Delta V(t)]\,dM\Big|\\
&\leqslant&\eta_9\cdot(||\Delta^2V(t)||_2^{\frac{1}{8}}+||D\Delta V(t)||_2^{\frac{1}{8}}+1)\\
&&\cdot(||\Delta^2V(t)||_2^{\frac{5}{8}}+||D\Delta V(t)||_2^{\frac{5}{8}}+1)\cdot||D\Delta V(t)||_2^{\frac{3}{8}}\cdot||\Delta^2V(t)||_2\\
&\lesssim&\eta_9\cdot(||\Delta^2V(t)||_2^{\frac{3}{4}}+||D\Delta V(t)||_2^{\frac{3}{4}}+1)\cdot||D\Delta V(t)||_2^{\frac{3}{8}}\cdot||\Delta^2V(t)||_2\\
&\leqslant&\eta_9\cdot||\Delta^2V(t)||_2^{\frac{7}{4}}\cdot||D\Delta V(t)||_2^{\frac{3}{8}}+\eta_9\cdot||\Delta^2V(t)||_2\cdot||D\Delta V(t)||_2^{\frac{9}{8}}\\
&&+\eta_9\cdot||\Delta^2V(t)||_2\cdot||D\Delta V(t)||_2^{\frac{3}{8}}\\
&\lesssim&\eta_9\cdot\Big(\varepsilon\cdot||\Delta^2V(t)||_2^2+\frac{1}{\varepsilon}\cdot||D\Delta V(t)||_2^3\Big)+\eta_9\cdot\Big(\varepsilon\cdot||\Delta^2V(t)||_2^2+\frac{1}{\varepsilon}\cdot||D\Delta V(t)||_2^{\frac{9}{4}}\Big)\\
&&+\eta_9\cdot\Big(\varepsilon\cdot||\Delta^2V(t)||_2^2+\frac{1}{\varepsilon}\cdot||D\Delta V(t)||_2^{\frac{3}{4}}\Big)\\
&\lesssim&\eta_9\cdot\Big(\varepsilon\cdot||\Delta^2V(t)||_2^2+\frac{1}{\varepsilon}\cdot||D\Delta V(t)||_2^4+\frac{1}{\varepsilon}\Big),
\end{eqnarray*}
where $\eta_9$ depends on $\mathcal{R}^M$, $\mathcal{R}^E$ and their covariant derivatives, $||V_0||_{\infty}, ||V_0||_2, ||DV_0||_2$, $||\Delta V_0||_2$ and $t$.

Moreover, there is a universal constant $\kappa_1$ such that
\begin{eqnarray*}
&&\Big|\int_M\langle\Delta\big[|V(t)|^2\cdot V(t)\big],\Delta^2 V(t)\rangle\,dM\Big|\\
&\leqslant&\frac{1}{4}||\Delta^2V(t)||_2^2+\kappa_1\cdot||V(t)||^4_{\infty}\cdot||\Delta V(t)||_2^2+\kappa_1\cdot||V(t)||^2_{\infty}\cdot||DV(t)||_4^4.
\end{eqnarray*}
Recalling (\ref{5.1}) and (\ref{56}) we obtain
\begin{eqnarray}\label{5.8}
&&\Big|\int_M\langle\Delta\big[|V(t)|^2\cdot V(t)\big],\Delta^2 V(t)\rangle\,dM\Big|\leqslant\frac{1}{4}||\Delta^2V(t)||_2^2+\kappa_1\cdot||V_0||^4_{\infty}\cdot\gamma_5\\
&&+\kappa_1\cdot||V_0||^2_{\infty}\cdot\eta_4\cdot(||D\Delta V(t)||_2+||\Delta V(t)||_2+||DV(t)||_2+||V(t)||_2)\cdot||DV(t)||_2^3.\nonumber
\end{eqnarray}
Substituting the upper bounds of $||\Delta V(t)||_2$, $||DV(t)||_2$ and $||V(t)||_2$ into (\ref{5.8}) gives
\begin{eqnarray*}
&&\Big|\int_M\langle\Delta\big[|V(t)|^2\cdot V(t)\big],\Delta^2 V(t)\rangle\,dM\Big|\\
&\leqslant&\frac{1}{4}||\Delta^2V(t)||_2^2+\kappa_2\cdot(||D\Delta V(t)||_2+1)\\
&\leqslant&\frac{1}{4}||\Delta^2V(t)||_2^2+\frac{\kappa_2}{2}\cdot(||D\Delta V(t)||^2_2+3),
\end{eqnarray*}
where $\kappa_2$ relies on $\mathcal{R}^M$, $\mathcal{R}^E$ and their covariant derivatives, $||V_0||_{\infty}, ||V_0||_2, ||DV_0||_2$, $||\Delta V_0||_2$ and $t$.

In conclusion,
\begin{eqnarray*}
&&\frac{1}{2}\frac{d}{dt}||D\Delta V(t)||_2^2+||\Delta^2V(t)||_2^2+\lambda\cdot||D\Delta V(t)||_2^2\\
&\leqslant&\eta'_9\cdot\Big(\varepsilon\cdot||\Delta^2V(t)||_2^2+\frac{1}{\varepsilon}\cdot||D\Delta V(t)||_2^4+\frac{1}{\varepsilon}\Big)+\frac{1}{4}||\Delta^2V(t)||_2^2+\frac{\kappa_2}{2}\cdot(||D\Delta V(t)||^2_2+3),
\end{eqnarray*}
where $\eta'_9$ depends on $\mathcal{R}^M$, $\mathcal{R}^E$ and their covariant derivatives, $||V_0||_{\infty}, ||V_0||_2, ||DV_0||_2$, $||\Delta V_0||_2$ and $t$. Let $\varepsilon$ be small enough. From Young's inequality it follows that
\begin{eqnarray*}
\frac{1}{2}\frac{d}{dt}||D\Delta V(t)||_2^2+\frac{1}{2}||\Delta^2V(t)||_2^2+\lambda\cdot||D\Delta V(t)||_2^2\leqslant\eta_{10}\cdot(1+||D\Delta V(t)||_2^4),
\end{eqnarray*}
where $\eta_{10}$ is dependent of $\eta_9$ and $\kappa_2$. Since of (\ref{Eq:1}), the generalized Gronwall's inequality implies
$$||D\Delta V(t)||_2^2\leqslant\gamma_6,$$
here $\gamma_6$ relies on $||D\Delta V_0||_2$, $\eta_{10}$, $\gamma_4$, $\gamma_5$ and $t$. Substituting the upper bounds of $||\Delta V(t)||^2_2$ and $||D\Delta V(t)||_2^2$ into (\ref{5.4}) and (\ref{5.5}) yields
$$||DV(t)||_{H^2}\leqslant\gamma_7$$
which implies
$$||DV(t)||_{\infty}\lesssim\gamma_7,$$
where $\gamma_7$ relies on $\mathcal{R}^M$, $\mathcal{R}^E$ and their covariant derivatives, $||V_0||_{\infty}, ||V_0||_2, ||DV_0||_2$, $||\Delta V_0||_2$, $||D\Delta V_0||_2$ and $t$.\\

Now we return to prove Theorem \ref{thm9}. Reviewing (\ref{37}), we know that the key is to estimate
$$\sum\limits_{i+j=k-1}a_{ij}\cdot||D^{k+1}V||_2\cdot\Big|\Big||D^{i+1}V|\cdot|D^{j+1}V|\Big|\Big|_2.$$
H\"older inequality yields
$$\Big|\Big||D^{i+1}V|\cdot|D^{j+1}V|\Big|\Big|_2\leqslant||D^{i+1}V||_{\frac{2k-2}{i}}\cdot||D^{j+1}V||_{\frac{2k-2}{k-1-i}}.$$
From Theorem \ref{thm6}, it follows that
$$||D^{i+1}V||_{\frac{2k-2}{i}}\lesssim||D^kV||_2^{\frac{i}{k-1}}\cdot||DV||_{\infty}^{\frac{k-1-i}{k-1}}$$
and
$$||D^{j+1}V||_{\frac{2k-2}{k-1-i}}\lesssim||D^kV||_2^{\frac{k-1-i}{k-1}}\cdot||DV||^{\frac{i}{k-1}}_{\infty}.$$
So one can get
\begin{eqnarray}\label{5.9}
\Big|\Big||D^{i+1}V|\cdot|D^{j+1}V|\Big|\Big|_2\lesssim||D^kV||_2\cdot||DV||_{\infty}.
\end{eqnarray}
Substituting (\ref{5.9}), (\ref{39}), (\ref{40}) and (\ref{41}) into (\ref{37}) we obtain
\begin{eqnarray*}
&&\frac{1}{2}\frac{d}{dt}||D^kV||_2^2+||D^{k+1}V||_2^2\\
&\leqslant&\kappa_3\cdot||D^{k+1}V||_2\cdot||D^kV||_2\cdot||DV||_{\infty}+G_k\cdot(1+||V_0||_{\infty}+||V_0||^2_{\infty})\cdot||V||^2_{H^k}\\
&\leqslant&\frac{1}{2}\cdot||D^{k+1}V||_2^2+\frac{\kappa_3^2}{2}\cdot||D^kV||_2^2\cdot||DV||^2_{\infty}+G_k\cdot(1+||V_0||_{\infty}+||V_0||^2_{\infty})\cdot||V||^2_{H^k},
\end{eqnarray*}
where $\kappa_3$ is a universal constant. Note the fact
$$||DV(t)||_{\infty}\lesssim\gamma_7\s\s\s\s\forall\s t\in[0,T^*).$$
Summing $k$ from 0 to $N$ and applying Gronwall's inequality we are led to
$$
||V(t)||^2_{H^N}+\int_0^t||V(s)||^2_{H^{N+1}}\,ds\leqslant C_N(||V_0||_{H^N},t,\gamma_7,||V_0||_{\infty},\kappa_3,G_0,\cdots,G_N).
$$
The remaining part of the proof of Theorem \ref{thm9} is as the same as that of Theorem \ref{thm5}. So we omit it. This completes the proof.\endproof

{}

\vspace{1.0cm}

Boling Guo

{\small\it Institute of Applied Physics and Computational Mathematics, China Academy of Engineering Physics, Beijing, 100088, P. R. China}

{\small\it Email: gbl@iapcm.ac.cn}\\

Zonglin Jia

{\small\it Institute of Applied Physics and Computational Mathematics, China Academy of Engineering Physics, Beijing, 100088, P. R. China}

{\small\it Email: 756693084@qq.com}\\


\begin{thebibliography}{2}
\bibitem{Baker}Charles Baker, The mean curvature flow of submanifolds of high codimension, Ph.D. thesis, Australian National University (2010); arXiv: 1104.4409.
\bibitem{Hamiltion} Richard S. Hamilton, Three-manifolds with positive Ricci curvature, J. Differential. Geometry, 17 (1982) 255-306.
\bibitem{DW1} Weiyue Ding, Youde Wang; \emph{Local Schr\"{o}dinger flow into K\"{a}hler manifolds}, Science in China, Series A, 2001, Vol. 44, No. 11.
\bibitem{GLZ} Boling Guo, Qiaoxin Li, Ming Zeng, Global smooth solutions of the Landau-Lifshitz-Bloch equation, preprint.
\bibitem{J} Zonglin Jia, Local strong solution to General Landau-Lifshitz-Bloch equation, arXiv:1802.00144.
\bibitem{L} Kim Ngan Le, Weak solutions of the Landau-Lifshitz Bloch equation, J. Differential Equations 261 (2016) 6699-6717
\end{thebibliography}
\end{document}